\documentclass[10pt]{article}
\usepackage{epsf}
\usepackage{amsmath}

\allowdisplaybreaks

\usepackage[showframe=false]{geometry}
\usepackage{changepage}

\usepackage{epsfig}
\usepackage{amssymb}

\usepackage{amsthm}
\usepackage{setspace}
\usepackage{cite}
\usepackage{mcite}

\usepackage{algorithmic}  
\usepackage{algorithm}

\usepackage{shadow}
\usepackage{fancybox}
\usepackage{fancyhdr}

\usepackage{color}
\usepackage[usenames,dvipsnames,svgnames,table]{xcolor}
\newcommand{\bl}[1]{\textcolor{blue}{#1}}

\definecolor{mypurple}{rgb}{.4,.0,.5}

\usepackage[hyphens]{url}

\usepackage[colorlinks=true,
            linkcolor=black,
            urlcolor=blue,
            citecolor=purple]{hyperref}

\usepackage{breakurl}

\def\w{{\bf w}}

\def\y{{\bf y}}

\def\x{{\bf x}}

\def\x{{\mathbf x}}

\def\w{{\bf w}}

\def\x{{\bf x}}
\def\y{{\bf y}}

\def\f{{\bf f}}

\def\diag{{\rm diag}\,}

\def\be{\begin{equation}}
\def\ee{\end{equation}}
\def\ba{\left[\begin{array}}
\def\ea{\end{array}\right]}

\def\w{{\bf w}}

\def\x{{\bf x}}
\def\y{{\bf y}}

\def\f{{\bf f}}

\def\1{{\bf 1}}
\def\oness{{\bf 1}^{(s)}}

\def\g{{\bf g}}
\def\0{{\bf 0}}

\def\erf{\mbox{erf}}
\def\erfc{\mbox{erfc}}
\def\erfi{\mbox{erfi}}
\def\diag{\mbox{diag}}

\def\mR{{\mathbb R}}
\def\mN{{\mathbb N}}
\def\mE{{\mathbb E}}

\def\phiint{\phi_{int}}
\def\phiext{\phi_{ext}}

\def\calF{{\cal F}}

\def\lp{\left (}
\def\rp{\right )}

\newtheorem{theorem}{Theorem}

\begin{document}

\begin{singlespace}

\title {Random linear systems with sparse solutions -- finite dimensions 
}
\author{
\textsc{Mihailo Stojnic
\footnote{e-mail: {\tt flatoyer@gmail.com}} }}
\date{}
\maketitle

\centerline{{\bf Abstract}} \vspace*{0.1in}

In our companion work \cite{Stojnicl1RegPosasymldp} we revisited random under-determined linear systems with sparse solutions. The main emphasis was on the performance analysis of the $\ell_1$ heuristic in the so-called asymptotic regime, i.e. in the regime where the systems' dimensions are large. Through an earlier sequence of work \cite{DonohoPol,DonohoUnsigned,StojnicCSetam09,StojnicUpper10}, it is now well known that in such a regime the $\ell_1$ exhibits the so-called \emph{phase transition} (PT) phenomenon. \cite{Stojnicl1RegPosasymldp} then went much further and established the so-called \emph{large deviations principle} (LDP) type of behavior that characterizes not only the breaking points of the $\ell_1$'s success but also the behavior in the entire so-called \emph{transition zone} around these points. Both of these concepts, the PTs and the LDPs, are in fact defined so that one can use them to characterize the asymptotic behavior. In this paper we complement the results of \cite{Stojnicl1RegPosasymldp} by providing an exact detailed analysis in the non-asymptotic regime. Of course, not only are the non-asymptotic results complementing those from \cite{Stojnicl1RegPosasymldp}, they actually are the ones that ultimately fully characterize the $\ell_1$'s behavior in the most general sense. We introduce several novel high-dimensional geometry type of strategies that enable us to eventually determine the $\ell_1$'s behavior. The main point will be the simplicity of the introduced techniques that then makes them fully applicable to a variety of other problems. In this introductory paper though we limit ourselves to the analysis of the $\ell_1$ to ensure the clarity/simplicity of the exposition and to respect the chronology of the creation of all of these results. Quite a few way more advanced concepts that fundamentally rely on what is presented here we will present separately in several companion papers.

\vspace*{0.25in} \noindent {\bf Index Terms: Finite dimensions; linear systems of equations; sparse solutions;
$\ell_1$-heuristic}.

\end{singlespace}

\section{Introduction}
\label{sec:back}

The problems that we will study below are directly related to what we studied in \cite{Stojnicl1RegPosasymldp}. Instead of repeating many of the introductory remarks from \cite{Stojnicl1RegPosasymldp} related to the description of the problems, their importance, and a large body of prior work that relates to it, we will often assume a decent level of familiarity with what was done in \cite{Stojnicl1RegPosasymldp} and in a bit lengthly line of work initiated in \cite{StojnicCSetamBlock09,StojnicUpperBlock10,StojnicICASSP09block,StojnicJSTSP09} and continued in many of our papers that followed afterwards. Consequently, we will here mostly focus on those explanations that we deem the most necessary to enable the smooth following of the main ideas and will often refer to \cite{Stojnicl1RegPosasymldp} for more thorough explanations. Along the same lines, below we will just briefly redescribe the linear system problems and then rather quickly switch to the presentation of the key topics of interest here.

As usual, we start by putting everything on the right mathematical track. Let $A$ be an $m\times n$ ($m\leq n$) dimensional real matrix (throughout the paper $A$ will on occasion be referred to as the system matrix) and let $\tilde{\x}$ be an $n$ dimensional real vector (for short we say, $A\in \mR^{m\times n}$ and $\tilde{\x}\in \mR^{n}$). $\tilde{\x}$ will be called $k$-sparse if no more than $k$ of its entries are different from zero. Let $\y$ be such that
\begin{equation}
\y\triangleq A\tilde{\x}. \label{eq:defy}
\end{equation}
The standard linear system problem is to determine $\tilde{\x}$ if $A$ and $\y$ in (\ref{eq:defy}) are given. This is then typically written as the following problem over unknown $\x$
\begin{equation}
A\x=\y. \label{eq:system}
\end{equation}
Algebraic properties of $A$ (dimensions and rank) determine if $\x$ can match $\tilde{\x}$. Here we will work with the full rank systems (i.e. we will assume that $A$ has full rank either deterministically or, when random, statistically). We will work in an under-determined regime ($m\leq n$) and $\tilde{\x}$ will be assumed to be $k$ sparse and the relation between $k$, $m$, and $n$ will be assumed such that the $k$-sparse $\x$ that solves (\ref{eq:system}) is unique and that there is no sparser $\x$ that solves (\ref{eq:system}) (again, either deterministically or, when $A$ is random, statistically). (\ref{eq:system}) is often also rewritten as the following optimization problem
\begin{eqnarray}
\mbox{min} & & \|\x\|_{0}\nonumber \\
\mbox{subject to} & & A\x=\y, \label{eq:l0gen}
\end{eqnarray}
where $\|\x\|_{0}$ is the so-called $\ell_0$ (quasi) norm of vector $\x$ ($\|\x\|_{0}$ will be viewed as the number of the nonzero entries of $\x$). Determining the sparsest $\x$ in (\ref{eq:l0gen}) (which we will technically call solving (\ref{eq:l0gen})) is of course typically considered as a not very easy task and finding provably good polynomial algorithms that would do so remains an extraordinary challenge (more on this can be found in e.g. \cite{StojnicCSetam09,Stojnicl1RegPosasymldp,StojnicReDirChall13}). Nonetheless, there are quite a few fairly successful algorithms developed over last several decades (see, e.g. \cite{JATGomp,NeVe07,DTDSomp,NT08,DaiMil08,DonMalMon09}) that could be utilized as solid heuristics. We view as mathematically the strongest the following $\ell_1$-optimization relaxation of (\ref{eq:l0gen})
\begin{eqnarray}
\mbox{min} & & \|\x\|_{1}\nonumber \\
\mbox{subject to} & & A\x=\y. \label{eq:l1}
\end{eqnarray}
Its popularity and importance are rooted in its polynomial complexity (it is a simple linear program) and excellent performance characteristics. An analysis of its performance characteristics is precisely what we will be interested in in the rest of this paper. Before we switch to showcasing the performance analysis that we developed we will just briefly mention several key milestones already achieved in the analysis of (\ref{eq:l1}).
The first breakthrough appeared in \cite{CRT,DOnoho06CS} where in a statistical scenario it was shown that (\ref{eq:l1}) can recover linear sparsity (recovering linear sparsity simply means that if the system's dimensions are large and $m$ is proportional to $n$ then there is a $k$ also proportional to $n$ such that the solutions of (\ref{eq:l0}) and (\ref{eq:l1}) coincide). These initial considerations were later on perfected in \cite{DonohoPol,DonohoUnsigned,StojnicCSetam09,StojnicUpper10,Stojnicl1RegPosasymldp}, where a qualitative description (such as linearity) was replaced with precise values of the system's dimensions/sparsity proportionalities. This effectively led to exact characterizations of the so-called $\ell_1$'s phase transitions in \cite{DonohoPol,DonohoUnsigned,StojnicCSetam09,StojnicUpper10} and to a complete characterization of a much stronger $\ell_1$'s large deviations principle in \cite{Stojnicl1RegPosasymldp}. While all these results were derived starting from finite dimensional systems they all were worked out so to emphasize the large dimensional asymptotic behavior. Here on the other hand, we will do exactly the opposite and attack the $\ell_1$'s performance analysis problem at its core considering the finite systems dimensions. To do so we will develop several novel high-dimensional geometry techniques. The main emphasis will be on their generality, simplicity, and elegance (these we found fairly useful and of great help in handling a variety of other problems).

We will organize the paper in a fairly simple way. We will split the presentation into two main sections. In the first one we will develop the main concepts by considering the positive $\ell_1$ variant of (\ref{eq:l1}) (see, e.g. \cite{Stojnicl1RegPosasymldp,StojnicCSetam09,StojnicUpper10}). In the second one we will then discuss all the needed upgrades so that the general (sometimes, we may also call it regular or standard) $\ell_1$ from (\ref{eq:l1}) can be handled as well.


\section{Positive $\ell_1$}
\label{sec:posl1}

As mentioned above, we will first consider a variant of (\ref{eq:l1}) that we will often call the positive $\ell_1$. Such a variant is typically used when the unknown $\x$ in (\ref{eq:l0}) (basically $\tilde{\x}$ in (\ref{eq:defy})) is structured a bit more beyond the standard sparsity. Namely, as above, one assumes that $\y$ was obtained through (\ref{eq:defy}) with $\tilde{\x}$ being not only $k$-sparse but also with components that are not negative (we will often refer to such vectors as nonnegative or positive). Assuming further that this is known beforehand and can be utilized in the algorithms' design, instead of (\ref{eq:l1}) one typically focuses on the following modification  (which, as mentioned above, we may often refer to as the \emph{nonnegative} or, a bit less precisely but possibly more naturally, the \emph{positive} $\ell_1$)
\begin{eqnarray}
\mbox{min} & & \|\x\|_{1}\nonumber \\
\mbox{subject to} & & A\x=\y\nonumber \\
&& \x\geq 0. \label{eq:l1nonn}
\end{eqnarray}
We will think of (\ref{eq:l1nonn}) as a relaxation of the following a priori known to have nonnegative sparse solution version of (\ref{eq:l0gen})
\begin{eqnarray}
\mbox{min} & & \|\x\|_{0}\nonumber \\
\mbox{subject to} & & A\x=\y\nonumber \\
&& \x\geq 0. \label{eq:l0}
\end{eqnarray}
As (\ref{eq:l1}), (\ref{eq:l1nonn}) is also a linear program and solvable in polynomial time. Moreover, as proven in \cite{DonohoPol,DonohoUnsigned,StojnicCSetam09,StojnicUpper10,Stojnicl1RegPosasymldp}, it has similar performance characteristics as the standard $\ell_1$ from (\ref{eq:l1}). Differently from putting an emphasis on the asymptotic/infinite dimensional setting, as was done in \cite{DonohoPol,DonohoUnsigned,StojnicCSetam09,StojnicUpper10,Stojnicl1RegPosasymldp}, below we will focus on finite dimensions and provide a performance analysis of (\ref{eq:l1nonn}) in such a scenario. As was the case in \cite{StojnicCSetam09,StojnicUpper10,Stojnicl1RegPosasymldp}, we start things off by recalling on a couple of results that we established in \cite{StojnicCSetam09,StojnicICASSP09} (these were, of course, the cornerstones of the success of the analysis done in \cite{StojnicCSetam09} as well as in a large sequence of our work that followed later on).

For the concreteness of the exposition and without loss of generality we will assume that the elements $\x_{k+1},\x_{k+2},\dots,\x_{n}$ of $\x$ are equal to zero and that the elements $\x_{1},\x_{2},\dots,\x_k$ are positive. The following is a positive adaptation of the result proven for the general $\ell_1$ in \cite{StojnicCSetam09,StojnicICASSP09} (see also \cite{YinZhang05nonneg}).
\begin{theorem}(\cite{StojnicCSetam09,StojnicICASSP09} Nonzero elements of $\x$ have fixed signs and location)
Assume that an $m\times n$ system matrix $A$ is given. Let $\x$
be a $k$ sparse vector. Also let $\x_{k+1}=\x_{k+2}=\dots=\x_{n}=0$. Let the signs of $\x_{1},\x_{2},\dots,\x_k$ be fixed, say all positive. Further, assume that $\y=A\x$ and that $\w$ is
a $n\times 1$ vector such that $\w_i\geq 0, k+1\leq i\leq n$. If
\begin{equation}
(\forall \w\in \textbf{R}^n | A\w=0) \quad  -\sum_{i=1}^{k} \w_i<\sum_{i=k+1}^{n}\w_{i},
\label{eq:posthmcond1}
\end{equation}
then the solutions of (\ref{eq:l0}) and (\ref{eq:l1nonn}) coincide. Moreover, if
\begin{equation}
(\exists \w\in \textbf{R}^n | A\w=0) \quad  -\sum_{i=1}^{k} \w_i\geq \sum_{i=k+1}^{n}\w_{i},
\label{eq:posthmcond2}
\end{equation}
then there is an $\x$ from the above set of $\x$'s with fixed location of nonzero elements such that the solution of (\ref{eq:l0}) is not the solution of (\ref{eq:l1nonn}).
\label{thm:posthmregposcond}
\end{theorem}
To facilitate the exposition we set
\begin{equation}
C^+_w\triangleq\{\w\in \mR^n| \quad -\sum_{i=1}^k \w_i\geq \sum_{i=k+1}^{n}\w_{i},\w_i\geq 0,k+1\leq i\leq n\}.\label{eq:posdefSw}
\end{equation}
Assuming that $A$ is random and that it has say i.i.d. standard normal components (alternatively one can consider $A$ that has the null-space uniformly distributed in the corresponding Grassmanian) the failing condition given in (\ref{eq:posthmcond2}) can be utilized for performance characterization of (\ref{eq:l1nonn}) in the following way
\begin{equation}
p^+_{err}(k,m,n)=P(\exists \w\in \textbf{R}^n | A\w=0, \w_i\geq 0,k+1\leq i\leq n, \quad \mbox{and} \quad -\sum_{i=1}^{k} \w_i\geq \sum_{i=k+1}^{n}\w_{i})
=P(G^{(sub)}\cap C^+_w\neq \emptyset).\label{eq:posanal1}
\end{equation}
Clearly, $p^+_{err}(k,m,n)$ is the probability that the solution of (\ref{eq:l1nonn}) is not the a priori known to be positive solution of (\ref{eq:l0}) and $G^{(sub)}$ is a uniformly randomly chosen subspace from $G_{n,n-m}$, the Grassmanian of all $(n-m)$-dimensional subspaces of $\mR^n$. Now, it is rather obvious that $C^+_w$ is a polyhedral cone (from this point on we will often assume a fair amount of familiarity with some of the well known concepts in high-dimensional geometry; more on them though can be found in e.g. \cite{BG,SantaloBookIntGeom76}). The following is a direct consequence of a remarkable result of Santalo \cite{Santalo} (see also, e.g. \cite{PMM,AS,BG})
\begin{equation}
P(G^{(sub)}\cap C^+_w\neq \emptyset)=2\sum_{l=m+2j+1,j\in \mN_0}^{n} \sum_{F^{(l,+)}\in \calF^{(l,+)}}\phiint(0,F^{(l,+)})\phiext(F^{(l,+)},C^+_w),\label{eq:posanal2}
\end{equation}
where $\calF^{(l,+)}$ is the set of all $l$-faces of $C^+_w$ and $\phiint(\cdot,\cdot)$ and $\phiext(\cdot,\cdot)$ are the so-called internal and external angles, respectively (more on the definitions and importance of $l$-faces, internal and external angles can be found in e.g. \cite{BG}). Connecting (\ref{eq:posanal1}) and (\ref{eq:posanal2}) we then also have
\begin{equation}
p^+_{err}(k,m,n)=P(G^{(sub)}\cap C^+_w\neq \emptyset)=2\sum_{l=m+2j+1,j\in \mN_0}^{n} \sum_{F^{(l,+)}\in \calF^{(l,+)}}\phiint(0,F^{(l,+)})\phiext(F^{(l,+)},C^+_w).\label{eq:posanal3}
\end{equation}
(\ref{eq:posanal3}) then seems to easily establish a way to characterize the performance of the positive $\ell_1$. However, things are not as simple as they seem and (\ref{eq:posanal3}) is just a conceptual solution to the problem of determining $p^+_{err}(k,m,n)$. To have (\ref{eq:posanal3}) be fully operational one would need to be able to handle the angles $\phiint(\cdot,\cdot)$ and $\phiext(\cdot,\cdot)$. This is typically a very hard problem and very rarely one encounters scenarios where the angles can be computed explicitly (needless to say that this is essentially exactly what makes determining $p^+_{err}(k,m,n)$ fairly hard). Below we present, in our view, a fairly elegant way to handle the angles. To ensure the clarity and easiness of the exposition, the presentation will be split into two parts, the first one that relates to the internal angles and the second one that deals with the external angles.

\subsection{Internal angles}
\label{sec:posintang}

In this section we analyze the internal angles appearing in (\ref{eq:posanal3}), i.e. $\phiint(0,F^{(l,+)})$. We will first distinguish between two types of $l$-faces of $C^+_w$ and will split the set of all $l$-faces $\calF^{(l,+)}$ into two sets, $\calF^{(l,+)}_1$ and $\calF^{(l,+)}_2$. To that end we set
\begin{eqnarray}
I_l & =  & \{1,2,\dots,k\}\nonumber \\
I_r & =  & \{k+1,k+2,\dots,n\},\label{eq:posintanal1}
\end{eqnarray}
and write
\begin{equation}
\calF^{(l,+)}_1 \triangleq\{\w\in \mR^n| \quad -\sum_{i=1}^k \w_i= \sum_{i=k+1}^{n}\w_{i},\w_{I_r}\geq 0,\w_{I^{(l,+)}_1}=0,I^{(l,+)}_1\subset I_r,|I^{(l,+)}_1|=n-l-1\}, l\in \{k-1,k,\dots,n-1\},\label{eq:posintanal2}
\end{equation}
and
\begin{equation}
\calF^{(l,+)}_2 \triangleq\{\w\in \mR^n| \quad -\sum_{i=1}^k \w_i\geq \sum_{i=k+1}^{n}\w_{i},\w_{I_r}\geq 0,\w_{I^{(l,+)}_2}=0,I^{(l,+)}_2\subset I_r,|I^{(l,+)}_2|=n-l\}, l\in \{k,k+1,\dots,n\},\label{eq:posintanal3}
\end{equation}
where for $Q\subset \{1,2,\dots,n\}$, $\w_{Q}$ is a vector that contains elements of $\w$ indexed by the elements of $Q$ (ordered for example according to the same way the indexes in $Q$ are ordered) and $\w_Q\geq 0$ means that each element of $\w_Q$ is larger than or equal to zero and analogously $\w_Q=0$ means that each element of $\w_Q$ is equal to zero. It is not that hard to see that the cardinalities of sets $\calF^{(l,+)}_1$ and $\calF^{(l,+)}_2$ are given by
\begin{equation}
c^{(l,+)}_1\triangleq |\calF^{(l,+)}_1|=\binom{n-k}{n-l-1}, l\in \{k-1,k,\dots,n-1\},\label{eq:posintanal4}
\end{equation}
and
\begin{equation}
c^{(l,+)}_2\triangleq |\calF^{(l,+)}_2|=\binom{n-k}{n-l}, l\in \{k,k+1,\dots,n\}.\label{eq:posintanal5}
\end{equation}
(\ref{eq:posanal3}) can then be rewritten in the following way
\begin{eqnarray}
p^+_{err}(k,m,n) & = & 2\sum_{l=m+2j+1,j\in \mN_0}^{n} ( \sum_{F^{(l,+)}_1\in \calF^{(l,+)}_1}\phiint(0,F^{(l,+)}_1)\phiext(F^{(l,+)}_1,C^+_w)\nonumber \\
& & + \sum_{F^{(l,+)}_2\in \calF^{(l,+)}_2}\phiint(0,F^{(l,+)}_2)\phiext(F^{(l,+)}_2,C^+_w))\nonumber \\
& = & 2 ( \sum_{l=m+2j+1,j\in \mN_0}^{n-1} c^{(l,+)}_1\phiint(0,F^{(l,+)}_1)\phiext(F^{(l,+)}_1,C^+_w) \nonumber \\
& & + \sum_{l=m+2j+1,j\in \mN_0}^{n} c^{(l,+)}_2 \phiint(0,F^{(l,+)}_2)\phiext(F^{(l,+)}_2,C^+_w)),\nonumber \\
\label{eq:posintanal6}
\end{eqnarray}
where due to symmetry $F^{(l,+)}_1$ and $F^{(l,+)}_2$ are basically any of the elements from sets $\calF^{(l,+)}_1$ and $\calF^{(l,+)}_2$, respectively. For the concreteness we choose $I^{(l,+)}_1=\{l+2,l+3,\dots,n\}$ and $I^{(l,+)}_2=\{l+1,l+2,\dots,n\}$ and consequently
\begin{equation}
F^{(l,+)}_1 =\{\w\in \mR^n| \quad -\sum_{i=1}^k \w_i= \sum_{i=k+1}^{n}\w_{i},\w_{I_r}\geq 0,\w_{I^{(l,+)}_1}=0\}, l\in \{k-1,k,\dots,n-1\},\label{eq:posintanal7}
\end{equation}
and
\begin{equation}
F^{(l,+)}_2 =\{\w\in \mR^n| \quad -\sum_{i=1}^k \w_i\geq \sum_{i=k+1}^{n}\w_{i},\w_{I_r}\geq 0,\w_{I^{(l,+)}_2}=0\}, l\in \{k,k+1,\dots,n\}.\label{eq:posintanal8}
\end{equation}
Below we handle separately $\phiint(0,F^{(l,+)}_1)$ and $\phiint(0,F^{(l,+)}_2)$.

\subsubsection{Handling $\phiint(0,F^{(l,+)}_1)$}
\label{sec:posint1ang}

To compute $\phiint(0,F^{(l,+)}_1)$ we will utilize a simple yet very powerful parametrization strategy that we will call ``Gaussian coordinates in an orthonormal basis". As the strategy lends itself to a repetitive use we will describe it separately.

\underline{\textbf{\emph{Gaussian coordinates in an orthonormal basis}}}

The strategy assumes two steps: 1) Finding an orthonormal basis in the subspace where the angle is being computed; 2) Expressing the content of the angle in terms of the Gaussian coordinates of the computed orthonormal basis.

1) For the orthonormal basis we will use the column vectors of the following matrix
\begin{equation}
B_{int,1}=\begin{bmatrix}
            \begin{bmatrix}
              B \\
              \0_{1\times (l-1)}
            \end{bmatrix}  & \begin{bmatrix}
              -\1_{l\times 1} \\
              l
            \end{bmatrix}\frac{1}{\sqrt{l^2+l}}\\
            \0_{(n-l-1)\times (l-1)} & \0_{(n-l-1)\times 1}
          \end{bmatrix},\label{eq:posint1anal1}
\end{equation}
where $B$ is an $l\times (l-1)$ orthonormal matrix such that $\1_{1\times l} B=\0_{(l-1)\times 1}$, and $\1$ and $\0$ are matrices of all ones or zeros, respectively of the sizes given in their indexes. By the above definition it is rather clear that $B_{int,1}^TB_{int,1}=I$. Moreover, $F^{(l,+)}_1$ is indeed in the subspace spanned by the columns of $B_{int,1}$ since its normal vector
$\f^{(l,1)}=\begin{bmatrix}-\1_{1\times (l+1)} & \0_{1\times (n-l-1)}\end{bmatrix}^T$ does satisfy $(\f^{(l,1)})^T B_{int,1}=\0_{1\times l}$.

2) Now that we have an orthonormal basis we will operate in it through the standard normal (i.e. Gaussian) coordinates. In other words every point in the subspace spanned by the columns of $B_{int,1}$ we will express in terms of $B_{int,1}\g$ where $\g$ will be assumed to have $l$ i.i.d standard normal components. However, not any $\g$ will be allowed but rather only those for which $B_{int,1}\g\in F^{(l,+)}_1$. Basically, the following set of $\g$'s, $G^{(l,+)}_1$, will be allowed
\begin{equation}
G^{(l,+)}_1 =\{\g\in \mR^l| \w=B_{int,1}\g \quad \mbox{and}\quad \w_i\geq 0,k+1\leq i\leq l+1\}, l\in \{k-1,k,\dots,n-1\}.\label{eq:posint1anal2}
\end{equation}
Utilizing (\ref{eq:posint1anal1}), (\ref{eq:posint1anal2}) can be rewritten as
\begin{equation}
G^{(l,+)}_1 =\{\g\in \mR^l| \w_{1:(l+1)}=\begin{bmatrix}
            \begin{bmatrix}
              B \\
              \0_{1\times (l-1)}
            \end{bmatrix}  & \begin{bmatrix}
              -\1_{l\times 1} \\
              l
            \end{bmatrix}\frac{1}{\sqrt{l^2+l}}
          \end{bmatrix}\g \quad \mbox{and}\quad \w_i\geq 0,k+1\leq i\leq l+1\}, l\in \{k-1,k,\dots,n-1\}.\label{eq:posint1anal3}
\end{equation}
Now we can write for $\phiint(0,F^{(l,+)}_1)$ (basically a definition of the internal angle)
\begin{equation}
\phiint(0,F^{(l,+)}_1)=\frac{1}{(2\pi)^{\frac{l}{2}}}\int_{G^{(l,+)}_1}e^{-\frac{\g^T\g}{2}}d\g.\label{eq:posint1anal4}
\end{equation}
We will also change variables in the following way (clearly relying on (\ref{eq:posint1anal3}))
\begin{equation}
\w_{1:l}=\begin{bmatrix}
              B  &  -\1_{l\times 1}\frac{1}{\sqrt{l^2+l}}
          \end{bmatrix}\g.
          \label{eq:posint1anal5}
\end{equation}
From (\ref{eq:posint1anal5}) we then obtain
\begin{equation}
\g=\begin{bmatrix}
              B^T  \\  -\1_{1\times l}\sqrt{\frac{l+1}{l}}
          \end{bmatrix}\w_{1:l}.
          \label{eq:posint1anal6}
\end{equation}
From (\ref{eq:posint1anal6}) we also have
\begin{eqnarray}
\g^T\g & = & \w_{1:l}^T\begin{bmatrix}
              B^T  \\  -\1_{1\times l}\sqrt{\frac{l+1}{l}}
          \end{bmatrix}^T\begin{bmatrix}
              B^T  \\  -\1_{1\times l}\sqrt{\frac{l+1}{l}}
          \end{bmatrix}\w_{1:l}\nonumber \\
          &= &\w_{1:l}^T\begin{bmatrix}
              B^T  \\  -\1_{1\times l}\sqrt{\frac{1}{l}}
          \end{bmatrix}^T\begin{bmatrix}
              B^T  \\  -\1_{1\times l}\sqrt{\frac{1}{l}}
          \end{bmatrix}\w_{1:l}+\w_{1:l}^T( -\1_{1\times l})^T
             (-\1_{1\times l})\w_{1:l}\nonumber \\
             & = & \w_{1:l}^T \w_{1:l}+\w_{1:l}^T( -\1_{1\times l})^T
             (-\1_{1\times l})\w_{1:l} .
          \label{eq:posint1anal7}
\end{eqnarray}
From (\ref{eq:posint1anal3}) and (\ref{eq:posint1anal5}) we have that $\w_{k+1:l}\geq 0$. (\ref{eq:posint1anal3}) forces $\g_l\geq 0$ which by (\ref{eq:posint1anal6}) implies that $(-\1_{1\times l})\w_{1:l}\geq 0$. For the Jacobian of the change of variables in (\ref{eq:posint1anal5}) one easily has
\begin{equation}
J=\frac{1}{\sqrt{\det\lp\begin{bmatrix}
              B  &  -\1_{l\times 1}\frac{1}{\sqrt{l^2+l}}
          \end{bmatrix}^T\begin{bmatrix}
              B  &  -\1_{l\times 1}\frac{1}{\sqrt{l^2+l}}
          \end{bmatrix}\rp}}=\sqrt{l+1}.\label{eq:posint1anal7a}
\end{equation}
(\ref{eq:posint1anal4}) can now be rewritten as
\begin{eqnarray}
\phiint(0,F^{(l,+)}_1) &  = & \frac{J}{(2\pi)^{\frac{l}{2}}}\int_{-\1_{1\times l}\w_{1:l}\geq 0,\w_{k+1:l}\geq 0}e^{-\frac{\w_{1:l}^T\w_{1:l}+\w_{1:l}^T( -\1_{1\times l})^T
             (-\1_{1\times l})\w_{1:l}}{2}}d\w_{1:l}\nonumber \\
&  = & \frac{J2^{l-k}}{2^{l-k}(2\pi)^{\frac{l}{2}}}\int_{-\1_{1\times l}\w_{1:l}\geq 0,\w_{k+1:l}\geq 0}e^{-\frac{\w_{1:l}^T\w_{1:l}+\w_{1:l}^T( -\1_{1\times l})^T
             (-\1_{1\times l})\w_{1:l}}{2}}d\w_{1:l}\nonumber \\
&  = & \frac{J}{2^{l-k}}\int_{x\geq 0}f_{x^+}(x)e^{-\frac{x^2}{2}}dx,\label{eq:posint1anal8}
\end{eqnarray}
and $f_{x^+}(\cdot)$ is the pdf of the random variable $x^+=-\1_{1\times l}\w_{1:l}$, where $\w_i,1 \leq i\leq k$, are  i.i.d. standard normals, and $\w_i, k+1 \leq i\leq l$, are  i.i.d. standard half-normals. To determine $f_{x^+}(\cdot)$, we will first compute the characteristic function of $x^+$. To that end we have
\begin{equation}
E e^{itx^+}=\frac{2^{l-k}}{(2\pi)^{\frac{l}{2}}}\int_{\w_{k+1:l}\geq 0}e^{-\frac{\w_{1:l}^T\w_{1:l}}{2}-it\1_{1\times l}\w_{1:l}}d\w_{1:l}
=e^{-\frac{kt^2}{2}}\lp1-i\erfi\lp\frac{t}{\sqrt{2}}\rp\rp^{l-k} e^{-\frac{(l-k)t^2}{2}},\label{eq:posint1anal9}
\end{equation}
and
\begin{equation}
f_{x^+}(x)=\frac{1}{2\pi}\int_{-\infty}^{\infty} e^{-\frac{kt^2}{2}}\lp1-i\erfi\lp\frac{t}{\sqrt{2}}\rp\rp^{l-k} e^{-\frac{(l-k)t^2}{2}}e^{-itx}dt,\label{eq:posint1anal10}
\end{equation}
where
\begin{equation}
\erfi(y)=-i\erf(iy)=\frac{2}{\sqrt{\pi}}\int_{0}^{y}e^{\frac{z^2}{2}}dz.\label{eq:posint1anal10a}
\end{equation}
Combining (\ref{eq:posint1anal8}) and (\ref{eq:posint1anal10}) we obtain
\begin{eqnarray}
\phiint(0,F^{(l,+)}_1) &  = &  \frac{J}{2^{l-k}}\int_{x\geq 0}f_{x^+}(x)e^{-\frac{x^2}{2}}dx\nonumber \\
& = & \frac{J}{2^{l-k}}\int_{x\geq 0}\frac{1}{2\pi}\int_{-\infty}^{\infty} e^{-\frac{kt^2}{2}}\lp1-i\erfi\lp\frac{t}{\sqrt{2}}\rp\rp^{l-k} e^{-\frac{(l-k)t^2}{2}}e^{-itx}dt e^{-\frac{x^2}{2}}dx\nonumber \\
& = & \frac{J}{2^{l-k+1}\sqrt{2\pi}}\int_{-\infty}^{\infty}\frac{2}{\sqrt{2\pi}}\int_{x\geq 0} e^{-\frac{kt^2}{2}}\lp1-i\erfi\lp\frac{t}{\sqrt{2}}\rp\rp^{l-k} e^{-\frac{(l-k)t^2}{2}}e^{-itx} e^{-\frac{x^2}{2}}dxdt\nonumber \\
& = & \frac{J}{2^{l-k+1}\sqrt{2\pi}}\int_{-\infty}^{\infty} e^{-\frac{kt^2}{2}}\lp1-i\erfi\lp\frac{t}{\sqrt{2}}\rp\rp^{l-k+1} e^{-\frac{(l-k+1)t^2}{2}}dt\nonumber \\
& = & \frac{\sqrt{l+1}}{2^{l-k+1}\sqrt{2\pi}}\int_{-\infty}^{\infty} \lp1-i\erfi\lp\frac{t}{\sqrt{2}}\rp\rp^{l-k+1} e^{-\frac{(l+1)t^2}{2}}dt.
\label{eq:posint1anal11}
\end{eqnarray}

\subsubsection{Handling $\phiint(0,F^{(l,+)}_2)$}
\label{sec:posint2ang}

Computing $\phiint(0,F^{(l,+)}_2)$ is a bit easier than computing $\phiint(0,F^{(l,+)}_1)$. There is no need to change variables and find a new orthonormal basis. Instead one can immediately write the following analogue to (\ref{eq:posint1anal8})
\begin{eqnarray}
\phiint(0,F^{(l,+)}_2) &  = & \frac{1}{(2\pi)^{\frac{l}{2}}}\int_{-\1_{1\times l}\w_{1:l}\geq 0,\w_{k+1:l}\geq 0}e^{-\frac{\w_{1:l}^T\w_{1:l}}{2}}d\w_{1:l}\nonumber \\
&  = & \frac{2^{l-k}}{2^{l-k}(2\pi)^{\frac{l}{2}}}\int_{-\1_{1\times l}\w_{1:l}\geq 0,\w_{k+1:l}\geq 0}e^{-\frac{\w_{1:l}^T\w_{1:l}}{2}}d\w_{1:l}\nonumber \\
&  = & \frac{1}{2^{l-k}}\int_{x\geq 0}f_{x^+}(x)dx.\label{eq:posint2anal8}
\end{eqnarray}
A combination of (\ref{eq:posint1anal10}) and (\ref{eq:posint2anal8}) then gives
\begin{eqnarray}
\phiint(0,F^{(l,+)}_2) &  = &  \frac{1}{2^{l-k}}\int_{x\geq 0}f_{x^+}(x)dx\nonumber \\
& = & \frac{1}{2^{l-k}}\int_{x\geq 0}\frac{1}{2\pi}\int_{-\infty}^{\infty} e^{-\frac{kt^2}{2}}\lp1-i\erfi\lp\frac{t}{\sqrt{2}}\rp\rp^{l-k} e^{-\frac{(l-k)t^2}{2}}e^{-itx}dt dx.
\label{eq:posint2anal11}
\end{eqnarray}
The above is sufficient to compute $\phiint(0,F^{(l,+)}_2)$. However, one can be even more explicit
\begin{eqnarray}
\phiint(0,F^{(l,+)}_2)
& = & \frac{1}{2^{l-k}}\int_{x\geq 0}\frac{1}{2\pi}\int_{-\infty}^{\infty} e^{-\frac{kt^2}{2}}\lp1-i\erfi\lp\frac{t}{\sqrt{2}}\rp\rp^{l-k} e^{-\frac{(l-k)t^2}{2}}e^{-itx}dt dx\nonumber \\
& = & \frac{1}{2^{l-k}}\frac{1}{2\pi}\lim_{x\rightarrow \infty}\lim_{\epsilon\rightarrow 0_+}(\int_{-\infty}^{-\epsilon} \lp1-i\erfi\lp\frac{t}{\sqrt{2}}\rp\rp^{l-k} e^{-\frac{lt^2}{2}}\frac{(1-e^{-itx})}{it}dt\nonumber \\
& & +\int_{\epsilon}^{\infty} \lp1-i\erfi\lp\frac{t}{\sqrt{2}}\rp\rp^{l-k} e^{-\frac{lt^2}{2}}\frac{(1-e^{-itx})}{it}dt).
\label{eq:posint2anal12}
\end{eqnarray}


\subsection{External angles}
\label{sec:posextang}

Similarly to what we did above for the internal angles, we will below compute separately the external angles of faces that belong to $\calF^{(l,+)}_1$ and $\calF^{(l,+)}_2$. In other words, we will separately handle $\phiext(F^{(l,+)}_1,C^+_w)$ and $\phiext(F^{(l,+)}_2,C^+_w)$.

\subsubsection{Handling $\phiext(F^{(l,+)}_1,C^+_w)$}
\label{sec:posext1ang}

By the definition of the external angles, they represent the content/fraction of the subspace taken by the positive hull of the outward normals to the hyperplanes that meet at the given face. For face $F^{(l,+)}_1$ we then have that the corresponding positive hull is given as
\begin{equation}\label{eq:posext1anal1}
phull^{(l,+)}_{ext,1}\triangleq  -pos(e_{l+2},e_{l+3},\dots,e_{n},-\1_{n\times 1}),
\end{equation}
where $e_i\in \mR^n$ and its only nonzero component is at the $i$-th location and is equal to one ($pos(\cdot,\cdot,\dots,\cdot)$ obviously stands for the positive hull of the vectors inside the parenthesis). To compute $\phiext(F^{(l,+)}_1,C^+_w)$ we will again rely on the ``Gaussian coordinates in an orthonormal basis" strategy that we presented above. That means that we first select an orthonormal basis in which we represent the content of $phull^{(l,+)}_{ext,1}$ (for the simplicity of writing we will actually work with $-phull^{(l,+)}_{ext,1}$; due to symmetry there is really no difference in the resulting angles). The basis that we will work in will be the columns of the following matrix
\begin{equation}
B^{(l,+)}_{ext,1}=\begin{bmatrix}
e_{l+2} &  e_{l+3} & \dots & e_{n} & \begin{bmatrix}
                                       -\1_{(l+1)\times 1} \\
                                       \0_{(n-l-1)\times 1}
                                     \end{bmatrix}\frac{1}{\sqrt{l+1}}
          \end{bmatrix}.\label{eq:posext1anal2}
\end{equation}
It is rather obvious that $pos(e_{l+2},e_{l+3},\dots,e_{n},-\1_{n\times 1})$ (or $-phull^{(l,+)}_{ext,1}$) is indeed located in the subspace spanned by the columns of $B^{(l,+)}_{ext,1}$. What is a bit more tricky is to accurately describe $pos(e_{l+2},e_{l+3},\dots,e_{n},-\1_{n\times 1})$ through the coordinates corresponding to the basis vectors from $B^{(l,+)}_{ext,1}$. To that end let these coordinates be denoted as earlier by a vector $\g$ (this time though $\g\in \mR^{n-l}$). What we want is a set of $\g$'s, say $G^{(l,+)}_{ext,1}$, such that
\begin{equation}
G^{(l,+)}_{ext,1} =\{\g\in \mR^{n-l}| B^{(l,+)}_{ext,1}\g\in pos(e_{l+2},e_{l+3},\dots,e_{n},-\1_{n\times 1})\triangleq-phull^{(l,+)}_{ext,1} \}.\label{eq:posext1anal3}
\end{equation}
Now let
\begin{equation}
D_{ext,1}=\begin{bmatrix}
e_{l+2} &  e_{l+3} & \dots & e_{n} & -\1_{n\times 1}\end{bmatrix},\label{eq:posext1anal4}
\end{equation}
and
\begin{eqnarray}
pos(e_{l+2},e_{l+3},\dots,e_{n},-\1_{n\times 1}) & = & \{D_{ext,1}\g^{(D)}|\g^{(D)}\geq 0,\g^{(D)}\in \mR^{n-l}\}\nonumber \\
& = & \{\begin{bmatrix}
e_{l+2} &  e_{l+3} & \dots & e_{n} & -\1_{n\times 1}\end{bmatrix}\g^{(D)}|\g^{(D)}\geq 0,\g^{(D)}\in \mR^{n-l}\}.\label{eq:posext1anal5}
\end{eqnarray}
From (\ref{eq:posext1anal3}) and (\ref{eq:posext1anal5}) we have
\begin{equation}
G^{(l,+)}_{ext,1} =\{\g\in \mR^{n-l}| \exists \g^{(D)}\in\mR^{n-l}, \g^{(D)}\geq 0, \quad \mbox{and} \quad B^{(l,+)}_{ext,1}\g= D_{ext,1}\g^{(D)} \}.\label{eq:posext1anal6}
\end{equation}
First $(l+1)$ equations in $B^{(l,+)}_{ext,1}\g= D_{ext,1}\g^{(D)}$ imply that
\begin{equation}
\g_{n-l}\frac{1}{\sqrt{l+1}}=\g^{(D)}_{n-l}\geq 0.\label{eq:posext1anal7}
\end{equation}
For $j\in\{2,3,\dots,n-l\}$, $(l+j)$-th equation in $B^{(l,+)}_{ext,1}\g= D_{ext,1}\g^{(D)}$ and (\ref{eq:posext1anal7}) also imply that
\begin{eqnarray}
& & \g_{j-1}=\g^{(D)}_{j-1}-\g^{(D)}_{n-l}=\g^{(D)}_{j-1}-\g_{n-l}\frac{1}{\sqrt{l+1}}\nonumber \\
\Longleftrightarrow & & \g_{j-1}+\g_{n-l}\frac{1}{\sqrt{l+1}}=\g^{(D)}_{j-1}\geq 0 \nonumber \\
\Longleftrightarrow & & \g_{j-1}\geq -\g_{n-l}\frac{1}{\sqrt{l+1}}
\label{eq:posext1anal8}
\end{eqnarray}
A combination of (\ref{eq:posext1anal6}), (\ref{eq:posext1anal7}), and (\ref{eq:posext1anal8}) then gives
\begin{equation}
G^{(l,+)}_{ext,1} =\{\g\in \mR^{n-l}| \g_{n-l}\geq 0, \g_{j-1}\geq -\g_{n-l}\frac{1}{\sqrt{l+1}}, j\in\{2,3,\dots,n-l\} \}.\label{eq:posext1anal9}
\end{equation}
Now that we expressed/parameterized the content of $-phull^{(l,+)}_{ext,1}$ through the coordinates $\g$ of the basis $B^{(l,+)}_{ext,1}$ we switch to the second part of the strategy. This part assumes working with the independent standard normal, i.e. Gaussian, coordinates. Working with them then automatically gives the fraction of the subspace taken by the set they define (as earlier, this is trivial due to the rotational invariance of the standard normal distribution). To that end we finally have for $\phiext(F^{(l,+)}_1,C^+_w)$
\begin{eqnarray}
\phiext(F^{(l,+)}_1,C^+_w) &  = & \frac{1}{(2\pi)^{\frac{n-l}{2}}}\int_{\g\in G^{(l,+)}_{ext,1}}e^{-\frac{\g^T\g}{2}}d\g \nonumber \\
& = & \frac{1}{(2\pi)^{\frac{l}{2}}}\int_{\g_{n-l}\geq 0}e^{-\frac{\g_{n-l}^2}{2}}\lp \prod_{j=2}^{n-l}\frac{1}{(2\pi)^{\frac{l}{2}}}\int_{\g_{j-1}\geq -\g_{n-l}\frac{1}{\sqrt{l+1}}}e^{-\frac{\g_{j-1}^2}{2}}d\g_{j-1}\rp d\g_{n-l}\nonumber \\
& = & \frac{1}{(2\pi)^{\frac{l}{2}}}\int_{\g_{n-l}\geq 0}e^{-\frac{\g_{n-l}^2}{2}}\lp \frac{1}{2}\erfc\lp\frac{-\g_{n-l}}{\sqrt{2}\sqrt{l+1}}\rp\rp^{n-l-1} d\g_{n-l}.
\label{eq:posext1anal10}
\end{eqnarray}

\subsubsection{Handling $\phiext(F^{(l,+)}_2,C^+_w)$}
\label{sec:posext2ang}

Computing $\phiext(F^{(l,+)}_2,C^+_w)$ is very simple compared to the above presented computing of $\phiext(F^{(l,+)}_1,C^+_w)$. The positive hull of the outward normals for face $F^{(l,+)}_2$ is given as
\begin{equation}\label{eq:posext2anal1}
phull^{(l,+)}_{ext,2}\triangleq  -pos(e_{l+1},e_{l+2},\dots,e_{n}).
\end{equation}
One then automatically has
\begin{equation}
\phiext(F^{(l,+)}_2,C^+_w)  =  \frac{1}{2^{n-l}}.
\label{eq:posext2anal2}
\end{equation}
All of what we presented above is then enough to determine $p^+_{err}(k,m,n)$. We summarize the obtained results in the following theorem.
\begin{theorem}(Exact \emph{nonnegative} $\ell_1$'s performance characterization -- finite dimensions)
Let $A$ be an $m\times n$ matrix in (\ref{eq:system})
with i.i.d. standard normal components (or, alternatively, with the null-space uniformly distributed in the Grassmanian). Let
the unknown $\x$ in (\ref{eq:system}) be $k$-sparse. Further, let all elements of $\x$ be nonnegative and let that be a priori known.
Let $p^+_{err}(k,m,n)$ be the probability that the solutions of (\ref{eq:l0}) and (\ref{eq:l1nonn}) do \emph{not} coincide. Then
\begin{eqnarray}
p^+_{err}(k,m,n) & = & 2 ( \sum_{l=m+2j+1,j\in \mN_0}^{n-1} c^{(l,+)}_1\phiint(0,F^{(l,+)}_1)\phiext(F^{(l,+)}_1,C^+_w)\nonumber \\
& & + \sum_{l=m+2j+1,j\in \mN_0}^{n} c^{(l,+)}_2 \phiint(0,F^{(l,+)}_2)\phiext(F^{(l,+)}_2,C^+_w)),
\label{eq:posfinalthm1}
\end{eqnarray}
where $c^{(l,+)}_1$, $ c^{(l,+)}_2$, $\phiint(0,F^{(l,+)}_1)$, $\phiint(0,F^{(l,+)}_2)$, $\phiext(F^{(l,+)}_1(C^+_w)$, and $\phiext(F^{(l,+)}_2(C^+_w)$ are as given in (\ref{eq:posintanal4}), (\ref{eq:posintanal5}), (\ref{eq:posint1anal11}), (\ref{eq:posint2anal12}), (\ref{eq:posext1anal10}), and  (\ref{eq:posext2anal2}), respectively.
\label{thm:posfinalperr}
\end{theorem}

\subsection{Simulations and theoretical results -- positive $\ell_1$}
\label{sec:posthnumresults}

\begin{figure}[htb]
\begin{minipage}[b]{.5\linewidth}
\centering
\centerline{\epsfig{figure=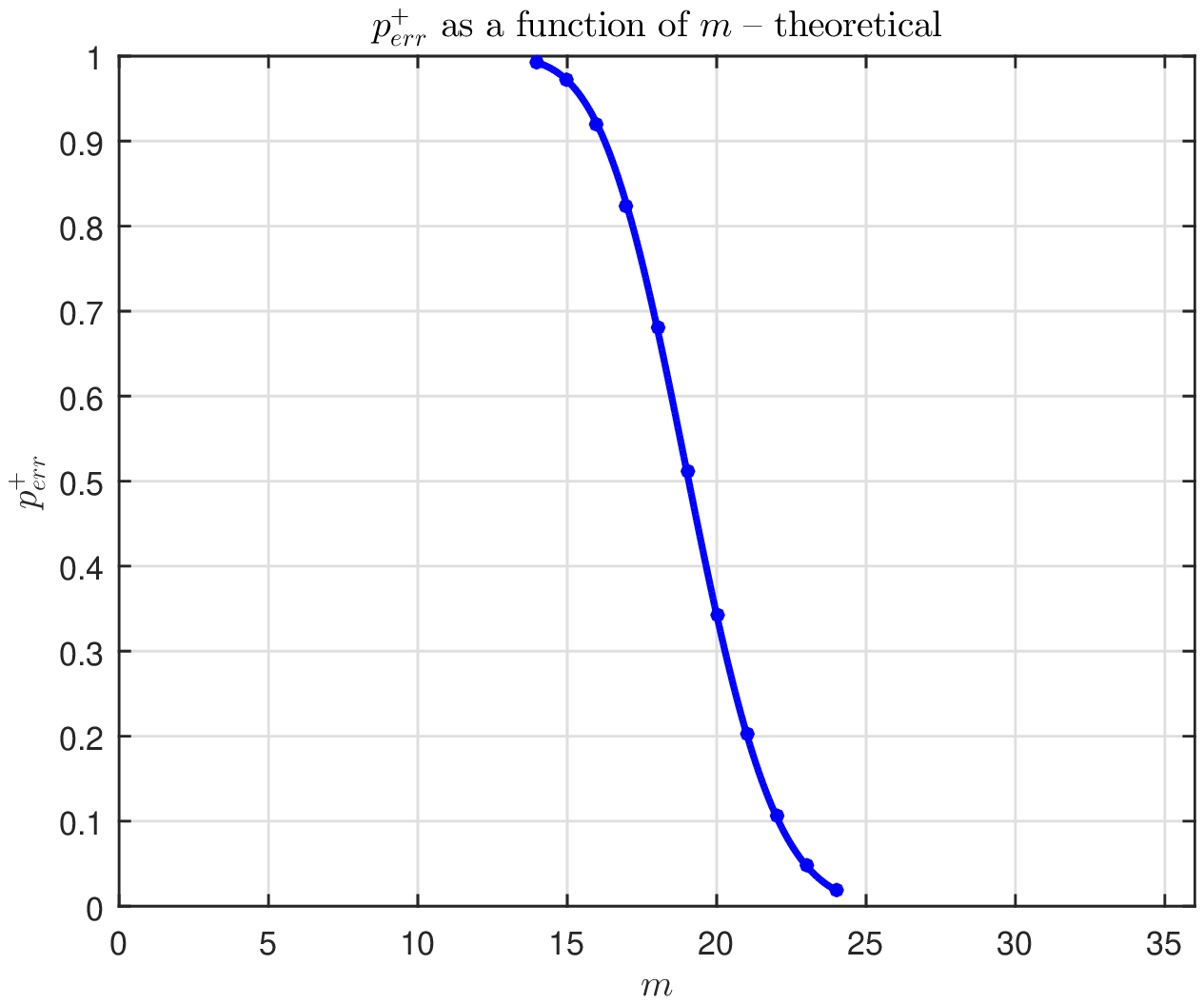,width=9cm,height=7cm}}
\end{minipage}
\begin{minipage}[b]{.5\linewidth}
\centering
\centerline{\epsfig{figure=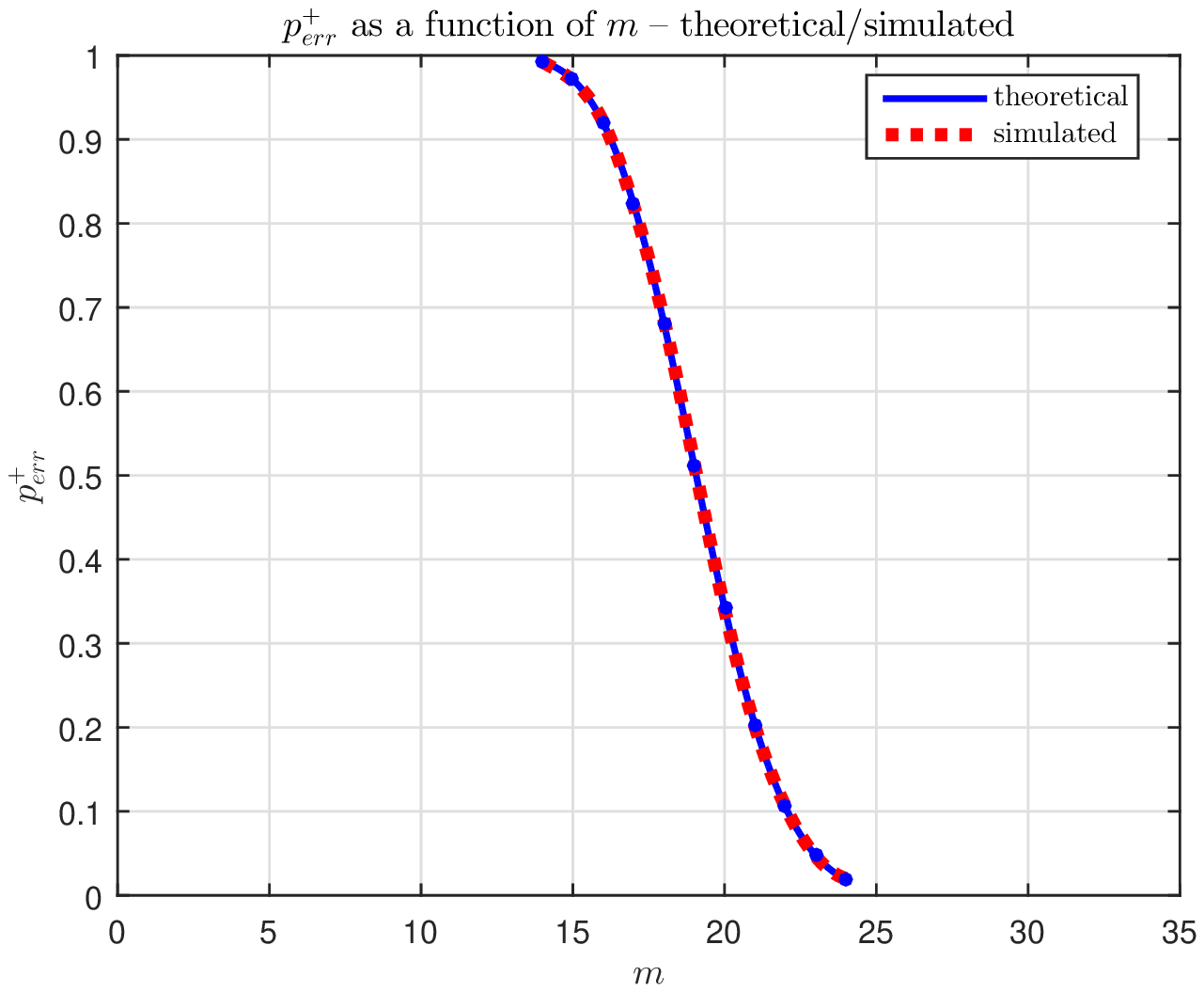,width=9cm,height=7cm}}
\end{minipage}
\caption{$p^+_{err}(k,m,n)$ as a function of $m$; left -- theory; right -- simulations}
\label{fig:l1regnonperr}
\end{figure}

In this section we will determine the concrete values for $p^+_{err}(k,m,n)$ based on what is proven in Theorem \ref{thm:posfinalperr}. In Figure \ref{fig:l1regnonperr} we show both, the simulated and the theoretical values for $p^+_{err}(k,m,n)$ (the theoretical values are, of course, obtained based on Theorem \ref{thm:posfinalperr}). We fixed $k=12$ and $n=36$ and varied/increased $m$ so that $p^+_{err}(k,m,n)$ changes from one to zero. Figure \ref{fig:l1regnonperr} is complemented by Table \ref{tab:l1regnonperrtab1} where we show the numerical values for $p^+_{err}(k,m,n)$ (again, both, simulated and theoretical) obtained for several concrete values of triplets $(k,m,n)$ (we also show the number of numerical experiments that were run as well as the number of them that did not result in having the solution of (\ref{eq:l1nonn}) match the a priori known to be nonnegative solution of (\ref{eq:l0})). We observe an excellent agreement between the simulated and theoretical results.

\begin{table}[h]
\caption{Simulated and theoretical results for $p^+_{err}(k,m,n)$; $k=12$, $n=36$}\vspace{.1in}
\hspace{-0in}\centering
\begin{tabular}{||c||c|c|c|c|c|c||}\hline\hline
$m$                      & $ 17 $ & $ 18 $ & $ 19 $ & $ 20 $ & $ 21 $ & $ 22 $ \\ \hline \hline
$\#$ of failures         & $ 10568 $ & $ 10944 $ & $ 9356 $ & $ 5574 $ & $ 2911 $ & $ 399 $ \\ \hline
$\#$ of repetitions      & $ 12869 $ & $ 15966 $ & $ 18196 $ & $ 16359 $ & $ 14464 $ & $ 3769 $ \\ \hline \hline
$p^+_{err}(k,m,n)$ -- simulation& $ \bl{\mathbf{0.8212}} $ & $ \bl{\mathbf{0.6855}} $ & $ \bl{\mathbf{0.5142}} $ & $ \bl{\mathbf{0.3407}} $ & $ \bl{\mathbf{0.2013}} $ & $ \bl{\mathbf{0.1059}} $ \\ \hline \hline
$p^+_{err}(k,m,n)$ -- theory    & $ \mathbf{0.8235} $ & $ \mathbf{0.6815} $ & $ \mathbf{0.5113} $ & $ \mathbf{0.3427} $ & $ \mathbf{0.2029} $ & $ \mathbf{0.1053} $ \\ \hline \hline
\end{tabular}
\label{tab:l1regnonperrtab1}
\end{table}

\subsection{Asymptotics}
\label{sec:posasym}

The above results are derived for concrete finite values of $k$, $m$, and $n$. As mentioned earlier, as such they establish the ultimate characterization of $p^+_{err}(k,m,n)$. When the systems dimensions grow larger (say in a linearly proportional fashion) the transition zone (where $p^+_{err}(k,m,n)$ changes from one to zero) gets relatively smaller and one ultimately, in an infinite dimensional setting, has the so-called phase-transition (PT) phenomenon. To properly introduce the PT phenomenon we will assume that $n$ is getting large and that $k=\beta n$ and $m=\alpha n$, where $\beta$ and $\alpha$ are fixed constants independent of $n$. Then the PT phenomenon, roughly speaking, boils down to the following: for a given $\alpha\in(0,1)$ there will be a $\beta^+_w\in(0,\alpha)$ such that for $\beta< \beta^+_w$ the a priori known to be nonnegative solution of (\ref{eq:l0}) is the solution of (\ref{eq:l1nonn}) with overwhelming probability. On the other hand, for $\beta\geq \beta^+_w$ there will be a $k$-sparse a priori known to be nonnegative $\x$ that is the solution of (\ref{eq:l0}) and is not the solution of (\ref{eq:l1nonn}) (as usual, under overwhelming probability we consider a probability that is no more than a number exponentially decaying in $n$ away from one). A full characterization of the PT phenomenon assumes finding $\beta^+_w$ for any $\alpha\in(0,1)$. A set of such points then determines what is often called the phase-transition curve (PT curve) in $(\alpha,\beta)$ plane (typically such curves are actually shown in $(\alpha,\beta/\alpha)$ plane so that one can emphasize the phase-transition dependence of the ratios $\frac{k}{m}$ and $\frac{m}{n}$).

The PT curve fully settles the PT phenomenon which as explained above deals with the determination of the so-called breaking points where (\ref{eq:l1nonn}) solves or does not solve (\ref{eq:l0}). A way more complete picture that describes the behavior of (\ref{eq:l1nonn})  in the entire transition zone around the breaking points can be obtained through the large deviations principle (LDP) type of analysis. This type of analysis instead of leaving the success/failure characterization of (\ref{eq:l1nonn}) in terms of overwhelming probabilities, insists on determining the exact rates of these probabilities' decaying parts. This rate is typically defined as $I^+_{ldp}(\alpha,\beta)=\lim_{n\rightarrow\infty}\frac{\log(p^+_{err}(k,m,n))}{n}$ and settling the LDP behavior of (\ref{eq:l1nonn}) assumes determining $I^+_{ldp}(\alpha,\beta)$ for any $(\alpha,\beta)\in(0,1)\times (0,1)$ for which the original system has unique solution.

The above analysis can be used to determine the LDP behavior and consequently the PT curve as well. Since settling the PT phenomenon is much easier than handling the entire LDP it can be done in a way that doesn't rely on the LDP results. In fact, the positive $\ell_1$'s PT behavior has been fully settled through the work of \cite{DT,DonohoSigned,StojnicCSetam09,StojnicUpper10}. The LDP on the other hand required a bit extra effort and we finally settled it in a companion paper \cite{Stojnicl1RegPosasymldp}. As is clear from \cite{Stojnicl1RegPosasymldp}, the LDP analysis is somewhat involved and a detailed discussion about it goes beyond what we study here which is the finite dimensional setup (exactly opposite of the asymptotic, infinite dimensional, one inherently assumed in the definitions of the PT and LDP phenomena). However, we below briefly sketch how one can bridge between the above analysis and what we will present in \cite{Stojnicl1RegPosasymldp}. In other words, we below show what problems a further analysis of the above $p^+_{err}(k,m,n)$ characterizations boils down to in an infinite dimensional setting (precisely these problems are, among other things, what we eventually solve in \cite{Stojnicl1RegPosasymldp}).

Until the end of this subsection we will assume the above mentioned so-called linear asymptotic regime, i.e. we will assume that $n$ is getting large and that $k=\beta n$ and $m=\alpha n$, where $\beta$ and $\alpha$ are fixed constants independent of $n$. Set $\rho\triangleq \lim_{n\rightarrow\infty}\frac{l}{n}$. Then, as $n\rightarrow\infty$, from (\ref{eq:posfinalthm1}) we have
\begin{equation}
\lim_{n\rightarrow \infty}\frac{\log(p^+_{err}(k,m,n))}{n}  =  \max\{\max_{\rho\geq \alpha}\lim_{n\rightarrow \infty} \frac{\log(\zeta^{(\infty,+)}_1)}{n},\max_{\rho\geq \alpha}\lim_{n\rightarrow \infty}\frac{\log(\zeta^{(\infty,+)}_2)}{n}\},
\label{eq:posasym1}
\end{equation}
where
\begin{eqnarray} \label{eq:posasym2}
  \lim_{n\rightarrow \infty} \frac{\log(\zeta^{(\infty,+)}_1)}{n} &=&
  \lim_{n\rightarrow \infty} \lp\frac{\log(c^{(l,+)}_{1})}{n}+
  \frac{\log(\phiint(0,F^{(l,+)}_1))}{n}+\frac{\log(\phiext(F^{(l,+)}_1,C^+_w))}{n}\rp \nonumber \\
  \lim_{n\rightarrow \infty} \frac{\log(\zeta^{(\infty,+)}_2)}{n} &=&
  \lim_{n\rightarrow \infty} \lp\frac{\log(c^{(l,+)}_{2})}{n}+
  \frac{\log(\phiint(0,F^{(l,+)}_2))}{n}+\frac{\log(\phiext(F^{(l,+)}_2,C^+_w))}{n}\rp.
\end{eqnarray}
From (\ref{eq:posintanal4}) and (\ref{eq:posintanal5}) we have
\begin{eqnarray}
\lim_{n\rightarrow \infty} \frac{\log(c^{(l,+)}_{1})}{n}=\lim_{n\rightarrow \infty}\frac{\log\binom{n-k}{n-l-1}}{n}=\lim_{n\rightarrow \infty}\frac{\log\binom{n-k}{n-l}}{n}=\lim_{n\rightarrow \infty} \frac{\log(c^{(l,+)}_{2})}{n}=-(1-\beta)H\lp\frac{1-\rho}{1-\beta}\rp,\label{eq:posasym3}
\end{eqnarray}
where
\begin{equation}\label{eq:posasym4}
  H(x)=x\log(x)+(1-x)\log(1-x).
\end{equation}
Now, from (\ref{eq:posint1anal8}) and (\ref{eq:posint2anal8}) we have
\begin{eqnarray}
\lim_{n\rightarrow \infty} \frac{\log(\phiint(0,F^{(l,+)}_1)}{n} &  = &
\lim_{n\rightarrow \infty} \frac{\log\lp\frac{\sqrt{l+1}}{(2\pi)^{\frac{l}{2}}}\int_{-\1_{1\times l}\w_{1:l}\geq 0,\w_{k+1:l}\geq 0}e^{-\frac{\w_{1:l}^T\w_{1:l}+\w_{1:l}^T( -\1_{1\times l})^T
             (-\1_{1\times l})\w_{1:l}}{2}}d\w_{1:l}\rp}{n}\nonumber \\
&  \leq  & \lim_{n\rightarrow \infty} \frac{\log\lp\frac{\sqrt{l+1}}{(2\pi)^{\frac{l}{2}}}\int_{-\1_{1\times l}\w_{1:l}\geq 0,\w_{k+1:l}\geq 0}e^{-\frac{\w_{1:l}^T\w_{1:l}}{2}}d\w_{1:l}\rp}{n}\nonumber \\
&  =  & \lim_{n\rightarrow \infty} \frac{\log(\sqrt{l+1})}{n}+\lim_{n\rightarrow \infty} \frac{\log\lp\frac{1}{(2\pi)^{\frac{l}{2}}}\int_{-\1_{1\times l}\w_{1:l}\geq 0,\w_{k+1:l}\geq 0}e^{-\frac{\w_{1:l}^T\w_{1:l}}{2}}d\w_{1:l}\rp}{n}\nonumber \\
& = & \lim_{n\rightarrow \infty} \frac{\log(\phiint(0,F^{(l,+)}_2)}{n}.\label{eq:posasym4a}
\end{eqnarray}
Now, let
\begin{equation}\label{eq:posasym4b}
  c_F(l)=\sqrt{\frac{\sqrt{l+1}-1}{l\sqrt{l+1}}},
\end{equation}
and
\begin{equation}\label{eq:posasym4c}
  B_F^{(l,+)}=\begin{bmatrix}\begin{bmatrix}
        I_{l\times l}-c_F(l)^2\1_{l\times 1}\1_{1\times l} & -\1_{l\times 1}\frac{1}{\sqrt{l+1}} & \\
        -\1_{1\times l}\frac{1}{\sqrt{l+1}} & -\frac{1}{\sqrt{l+1}}
        \end{bmatrix} & \0_{(l+1)\times (n-l-1)} \\
        \0_{(n-l-1)\times (l+1)} & \0_{(n-l-1)\times (n-l-1)}
      \end{bmatrix}.
\end{equation}
It is not that hard to check that $(B_F^{(l,+)})^TB_F^{(l,+)}=I_{n\times n}$ which implies
\begin{equation}\label{eq:posasym4ca}
\phiint(0,F^{(l,+)}_2)=\phiint(0,B_F^{(l,+)} F^{(l,+)}_2).
\end{equation}
We also recall that
\begin{eqnarray}\label{eq:posasym4d}
F^{(l,+)}_1 & = & \{\w\in \mR^n|-\sum_{i=1}^k \w_i = \sum_{i=k+1}^{l+1}\w_{i},\w_{k+1:l+1}\geq 0,\w_{l+2:n}= 0\} \nonumber \\
F^{(l,+)}_2 & = & \{\w\in \mR^n|-\sum_{i=1}^k \w_i\geq \sum_{i=k+1}^{l}\w_{i},\w_{k+1:l}\geq 0,\w_{l+1:n}= 0\}.
\end{eqnarray}
Let $\w^{(l,2)}\in F^{(l,+)}_2$. Then
\begin{equation} \label{eq:posasym4e}
\w^{(l,1)}=B_F^{(l,+)}\w^{(l,2)}=\begin{bmatrix}
        \w^{(l,2)}_{1:l}-c_F(l)^2\1_{l\times 1}\1_{1\times l}\w^{(l,2)}_{1:l}  \\
        -\1_{1\times l}\frac{1}{\sqrt{l+1}}\w^{(l,2)}_{1:l}  \\
        \0_{(n-l-1)\times 1}
      \end{bmatrix}.
\end{equation}
Clearly, $\w^{(l,1)}_{l+2:n}=\0_{(n-l-1)\times 1}$. By (\ref{eq:posasym4d}) we have $-\1_{1\times l}\frac{1}{\sqrt{l+1}}\w^{(l,2)}_{1:l}\geq 0$. Together with $\w^{(l,2)}_{k+1:l}\geq 0$ this implies that $\w^{(l,1)}_{k+1:l}\geq 0$. We also have
\begin{eqnarray} \label{eq:posasym4f}
-\1_{1\times n} \w^{(l,1)} & = & -\1_{1\times n} B_F^{(l,+)}\w^{(l,2)}\nonumber \\
& = & -\1_{1\times l}\w^{(l,2)}_{1:l}+c_F(l)^2\1_{1\times l}\1_{l\times 1}\1_{1\times l}\w^{(l,2)}_{1:l}
        +\1_{1\times l}\frac{1}{\sqrt{l+1}}\w^{(l,2)}_{1:l} \nonumber \\
& = & \1_{1\times l}\w^{(l,2)}_{1:l}(-1+c_F(l)^2\1_{1\times l}\1_{l\times 1}
       +\frac{1}{\sqrt{l+1}}) \nonumber \\
& = & \1_{1\times l}\w^{(l,2)}_{1:l}(-1+\frac{\sqrt{l+1}-1}{l\sqrt{l+1}}\1_{1\times l}\1_{l\times 1}
       +\frac{1}{\sqrt{l+1}}) \nonumber \\
& = & \1_{1\times l}\w^{(l,2)}_{1:l}(-1+\frac{\sqrt{l+1}-1}{\sqrt{l+1}}
       +\frac{1}{\sqrt{l+1}}) \nonumber \\
& = & 0.
\end{eqnarray}
By (\ref{eq:posasym4d}) this then implies that $\w^{(l,1)}=B_F^{(l,+)}\w^{(l,2)}\in F^{(l,+)}_1$ and ultimately
\begin{equation}\label{eq:posasym4g}
   B_F^{(l,+)}F^{(l,+)}_2 \in F^{(l,+)}_1.
\end{equation}
A combination of (\ref{eq:posasym4ca}) and (\ref{eq:posasym4g}) gives
\begin{equation}\label{eq:posasym4h}
\phiint(0,F^{(l,+)}_2)=\phiint(0,B_F^{(l,+)} F^{(l,+)}_2) \leq \phiint(0,F^{(l,+)}_1).
\end{equation}
Finally, combining (\ref{eq:posasym4a}) and (\ref{eq:posasym4h}) we obtain
\begin{eqnarray}
\lim_{n\rightarrow \infty} \frac{\log(\phiint(0,F^{(l,+)}_1)}{n}  =  \lim_{n\rightarrow \infty} \frac{\log(\phiint(0,F^{(l,+)}_2)}{n}.
\label{eq:posasym4i}
\end{eqnarray}
From (\ref{eq:posint2anal8}) we also have
\begin{eqnarray}
\phiint(0,F^{(l,+)}_2) &  = &  \frac{2^{l-k}}{2^{l-k}(2\pi)^{\frac{l}{2}}}\int_{-\1_{1\times l}\w_{1:l}\geq 0,\w_{k+1:l}\geq 0}e^{-\frac{\w_{1:l}^T\w_{1:l}}{2}}d\w_{1:l}\nonumber \\
&  = & \frac{1}{2^{l-k}}P(-\1_{1\times l}\w_{1:l}\geq 0),\label{eq:posasym6}
\end{eqnarray}
where on the right side of the last equality one can think of the elements of $\w_{1:k}$ as being the i.i.d. standard normals and the elements of $\w_{k+1:l}$ as being the i.i.d. standard half normals. By the definition of the large deviations principle we further have
\begin{eqnarray}
\lim_{n\rightarrow \infty} \frac{\log(\phiint(0,F^{(l,+)}_2))}{n}
&  = & \lim_{n\rightarrow \infty} \frac{\log(\frac{1}{2^{l-k}}P(-\1_{1\times l}\w_{1:l}\geq 0))}{n}\nonumber \\
& = &  \min_{\mu_y\geq 0} \lim_{n\rightarrow \infty} \frac{\log(\mE e^{-\mu_y\1_{1\times l}\w_{1:l}})}{n}-(\rho-\beta)\log(2)\nonumber \\
& = &  \min_{\mu_y\geq 0} \lp(\rho-\beta) \log\lp\mE e^{-\mu_y\w_{k+1}}\rp+\beta\frac{\mu_y^2}{2}\rp-(\rho-\beta)\log(2) \nonumber \\
& = &  \min_{\mu_y\geq 0} \lp(\rho-\beta) \log\lp\frac{2}{\sqrt{2\pi}}\int_{0}^{\infty} e^{-\frac{\w_{k+1}^2}{2}-\mu_y\w_{k+1}}d\w_{k+1}\rp+\beta\frac{\mu_y^2}{2}\rp-(\rho-\beta)\log(2) \nonumber \\
& = &  \min_{\mu_y\geq 0} \lp(\rho-\beta) \log\lp\erfc\lp\frac{\mu_y}{\sqrt{2}}\rp\rp+\rho\frac{\mu_y^2}{2}\rp-(\rho-\beta)\log(2) \nonumber \\
& = &  \min_{\mu_y\geq 0} \lp(\rho-\beta) \log(\erfc(\mu_y))+\rho\mu_y^2\rp-(\rho-\beta)\log(2).\label{eq:posasym7}
\end{eqnarray}
From (\ref{eq:posext1anal10}) we obtain
\begin{eqnarray}
\lim_{n\rightarrow \infty} \frac{\log(\phiext(F^{(l,+)}_1,C^+_w))}{n}
& = & \lim_{n\rightarrow \infty} \frac{\log\lp\frac{1}{(2\pi)^{\frac{l}{2}}}\int_{\g_{n-l}\geq 0}e^{-\frac{\g_{n-l}^2}{2}}\lp \frac{1}{2}\erfc\lp\frac{-\g_{n-l}}{\sqrt{2}\sqrt{l+1}}\rp\rp^{n-l-1} d\g_{n-l}\rp}{n} \nonumber \\
& = & \max_{\g_{n-l}\geq 0} \lp -\frac{\g_{n-l}^2}{2}+(1-\rho)\log\lp \frac{1}{2}\erfc\lp\frac{-\g_{n-l}}{\sqrt{2}\sqrt{\rho}}\rp\rp\rp \nonumber \\
& = & \max_{\g_{n-l}\geq 0} \lp -\rho \g_{n-l}^2+(1-\rho)\log\lp \frac{1}{2}\erfc(-\g_{n-l})\rp\rp,
\label{eq:posasym8}
\end{eqnarray}
and from (\ref{eq:posext2anal2}) we have
\begin{eqnarray}
\lim_{n\rightarrow \infty} \frac{\log(\phiext(F^{(l,+)}_2,C^+_w))}{n}
 =  \lim_{n\rightarrow \infty} \frac{\log\lp\frac{1}{2^{n-l}}\rp}{n}=-(1-\rho)\log(2).
\label{eq:posasym9}
\end{eqnarray}
For $\g_{n-l}=0$ in (\ref{eq:posasym8}) we have $\lim_{n\rightarrow \infty} \frac{\log(\phiext(F^{(l,+)}_1,C^+_w))}{n}=-(1-\rho)\log(2)$ which implies that  \begin{eqnarray}
\lim_{n\rightarrow \infty} \frac{\log(\phiext(F^{(l,+)}_2,C^+_w))}{n}
 \leq  \lim_{n\rightarrow \infty} \frac{\log(\phiext(F^{(l,+)}_1,C^+_w))}{n}.
\label{eq:posasym10}
\end{eqnarray}
A combination of (\ref{eq:posasym1}), (\ref{eq:posasym2}), (\ref{eq:posasym3}), (\ref{eq:posasym4i}), and (\ref{eq:posasym10}) gives
\begin{equation}
\lim_{n\rightarrow \infty}\frac{\log(p^+_{err}(k,m,n))}{n}  =  \max\{\max_{\rho\geq \alpha}\lim_{n\rightarrow \infty} \frac{\log(\zeta^{(\infty,+)}_1)}{n},\max_{\rho\geq \alpha}\lim_{n\rightarrow \infty}\frac{\log(\zeta^{(\infty,+)}_2)}{n}\}
=\max_{\rho\geq \alpha}\lim_{n\rightarrow \infty} \frac{\log(\zeta^{(\infty,+)}_1)}{n}.
\label{eq:posasym11}
\end{equation}
Relying on (\ref{eq:posasym2}), (\ref{eq:posasym3}), (\ref{eq:posasym7}), and (\ref{eq:posasym8}), one can rewrite (\ref{eq:posasym11}) in the following way
\begin{eqnarray}
\lim_{n\rightarrow \infty}\frac{\log(p^+_{err}(k,m,n))}{n}
 & = & \max_{\rho\geq \alpha}\lim_{n\rightarrow \infty} \frac{\log(\zeta^{(\infty,+)}_1)}{n}\nonumber \\
 & = & \max_{\rho\geq \alpha} (-(1-\beta)H\lp\frac{1-\rho}{1-\beta}\rp \nonumber \\
 & & + \min_{\mu_y\geq 0} \lp(\rho-\beta) \log(\erfc(\mu_y))+\rho\mu_y^2\rp-(\rho-\beta)\log(2)\nonumber \\.
& & + \max_{\g_{n-l}\geq 0} \lp -\rho \g_{n-l}^2+(1-\rho)\log\lp \frac{1}{2}\erfc(-\g_{n-l})\rp\rp).
\label{eq:posasym12}
\end{eqnarray}
For a given $\beta$, let $\alpha^+_w$ be the $\alpha$ that gives $\lim_{n\rightarrow \infty}\frac{\log(p^+_{err}(k,m,n))}{n} =0$ (it is rather trivial that such an $\alpha$ always exists; nonetheless, for the completeness, we in \cite{Stojnicl1RegPosasymldp} provide a straightforward rigorous proof of this fact). This is in fact how one can determine the phase transition values. Now, if $\alpha\geq \alpha_w$ the optimal $\rho$ in (\ref{eq:posasym12}) is trivially equal to $\alpha$. For example, if it wasn't, there would be a different $\alpha$, say $\alpha_2$, such that $\alpha_2=\rho>\alpha$ for which one would have that its $p^+_{err}(k,m,n)$ would be larger than the one that corresponds to $\alpha$. This of course is not possible since when the number of equations increases the probability of error cannot increase. On the other hand, if $\alpha\leq \alpha_w$ the optimal $\rho$ will again trivially be equal to $\alpha$. The reasoning is the same as above, the only difference is that now one looks at the complementary version of (\ref{eq:posanal2})
\begin{equation}
P(G^{(sub)}\cap C^+_w\neq \emptyset)=1-2\sum_{l=m-2j-1,j\in \mN_0,l\geq k-1} \sum_{F^{(l,+)}\in \calF^{(l,+)}}\phiint(0,F^{(l,+)})\phiext(F^{(l,+)},C^+_w)=1-p^+_{cor},\label{eq:posasym13}
\end{equation}
where $p^+_{cor}$ is the probability of being correct, i.e. the probability that the solution of (\ref{eq:l1nonn}) is the a priori known to be nonnegative solution of (\ref{eq:l0}) and its decay rate is given by (\ref{eq:posasym12}) with $\rho\leq \alpha$.
As above, if optimal $\rho$ is smaller than $\alpha\leq \alpha_w$, there would be a different $\alpha$, say $\alpha_2$, such that $\alpha_2=\rho<\alpha$ for which one would have that its $p^+_{cor}$ would be larger than the one that corresponds to $\alpha$. This again is not possible since when the number of equations decreases the probability of being correct (i.e. the probability that the solution of (\ref{eq:l1nonn}) is the a priori known to be nonnegative solution of (\ref{eq:l0})) cannot increase. (\ref{eq:posasym12}) is sufficient to fully determine numerically PT and LDP curves of the positive $\ell_1$. \cite{Stojnicl1RegPosasymldp}, however, goes way beyond that and determines explicit solutions to (\ref{eq:posasym12}).

\subsection{Face counting and projection survival}
\label{sec:posfacecountprsur}

In \cite{DonohoUnsigned} Donoho connected the success of (\ref{eq:l1}) in finding the solution of (\ref{eq:l0}) to the problem of counting faces of neighbourly polytopes. In \cite{DonohoSigned} such a strategy was adapted to the positive $\ell_1$ from (\ref{eq:l1nonn}). Roughly speaking, the reasoning in \cite{DonohoUnsigned} goes as follows: if a given $(k-1)$-face of the regular standard $n$-dimensional crosspolytope $C^{(n)}$ (that essentially corresponds to  a set of scaled $k$-sparse $\x$ in (\ref{eq:system})) so to say ``survives" the projection $A$, i.e. remains a face of the projected crosspolytope $C^{(n)}$, $AC^{(n)}$, then the solution of (\ref{eq:l1}) is the $k$-sparse solution of (\ref{eq:l0}) (which is located exactly on the given $(k-1)$-face of $C^{(n)}$). Assuming that $A$ is the orthoprojector on an $m$-dimensional subspace (uniformly chosen among all $m$-dimensional subspaces  of $\mR^n$)  \cite{DonohoPol} then in an asymptotic regime determined the precise values for $\alpha$ and $\beta$ so that an overwhelming number of $(k-1)$-faces of $C^{(n)}$ survives the projection. As mentioned above, the strategy was then adapted in \cite{DonohoSigned} so that one can eventually through \cite{DT} determine the positive $\ell_1$ PT curve. Instead of the standard corssolytope the standard simplex in $\mR^n$, $T^{(n-1)}$, was considered. This time though, one does not have an exact correspondence between the survival of a $(k-1)$-face of $T^{(n-1)}$ and the ability of the positive $\ell_1$ from (\ref{eq:l1nonn}) to handle the a priori known to have nonnegative solution (\ref{eq:l0}). Nonetheless, through considerations of a solid simplex $T^{(n-1)}_0$ (basically a convex hull of $0$ and $T^{(n-1)}$) and its an outward neighborliness \cite{DT,DonohoSigned} showed that survival of an overwhelming fraction of $(k-1)$-faces after projecting $T^{(n)}$ by $A$ implies the so-called weak equivalence of (\ref{eq:l0}) and (\ref{eq:l1nonn}). That then allowed the use of the known results \cite{AS,BetkeHenk,Ruben,BorockyHenk} related to the counting of the faces of the projected standard simplexes to determine in an asymptotic regime the positive $\ell_1$'s PT curve.

As we deal here with the finite dimensions, an exact correspondence would be welcome. Although, as mentioned above, this type of the correspondence in the sense of \cite{DonohoUnsigned} is lacking when one views the equivalence of (\ref{eq:l0}) and (\ref{eq:l1nonn}) through the random projection of the standard simplex, we actually found a way how the known results for counting faces of randomly projected simplexes from \cite{AS,BetkeHenk,Ruben,BorockyHenk} (basically the ones that were utilized in \cite{DT} for an asymptotic analysis and determination of the positive $\ell_1$'s phase transition) can still be of use even in the finite dimension scenarios. Namely, instead of looking at the performance of the positive $\ell_1$ one can look at its an alternative version
\begin{eqnarray}
\mbox{min} & & \|\x\|_{1}\nonumber \\
\mbox{subject to} & & A\x=\y\nonumber \\
& & \sum_{i=1}^{n}\x_i=\sum_{i=1}^{n}\tilde{\x}_i-\epsilon \nonumber \\
&& \x\geq 0, \label{eq:l1nonnsimp}
\end{eqnarray}
where $\epsilon>0$ (solving this problem would be possible of course if one is allowed a bit of a feedback in the form of available knowledge of the sum of the elements of the unknown vector $\x$). Adding the constraint $\sum_{i=1}^{n}\x_i=\sum_{i=1}^{n}\tilde{\x}_i-\epsilon$ essentially forces the problem to reside in $T^{(n-1)}$ and one now easily through a simple convex combination analysis has the one-one correspondence between projected faces survival and the equivalence of (\ref{eq:l0}) and (\ref{eq:l1nonn}).

Counting faces of the projected standard simplexes goes through the mechanism developed in \cite{AS}, which critically relies on the formulas similar to the ones we used above that are due to Grunbaum, McMullen, and ultimately Santalo. Here though, we don't need to count faces; instead we are looking for the probability that a $(k-1)$-face of $T^{(n-1)}$ fails to survive the projection by $A$ (this will give the probability of error $p^{+,s}_{err}(k,m,n)$, which is the probability that (\ref{eq:l1nonnsimp}) fails to uniquely produce the a priori known to be nonnegative solution of (\ref{eq:l0}) with a given sum $\sum_{i=1}^{n}\tilde{\x}_i$). To that end we write (following the definition of the Grunbaum's Grassman angle (see, e.g. \cite{Grunbaum68}) and its a characterization through \cite{Santalo} and the angle-sum relation of \cite{PMM})
\begin{equation}\label{eq:posfacecount1}
p^{+,s}_{err}(k,m,n)  =2\sum_{l=m+2j+1,j\in \mN_0}^{n-1} \sum_{F^{(l,+)}\in \calF^{(l,+)},F^{(k-1)}\in F^{(l,+)}}\phiint(F^{(k-1)},F^{(l,+)})\phiext(F^{(l,+)},T^{(n-1)}).
\end{equation}
As is trivially known, all $l$-faces of $T^{(n-1)}$ are also standard simplexes. Also, the number of $l$-faces that contain the given $(k-1)$-face is equal to $c^{(l,+)}_1$ and one can rewrite (\ref{eq:posfacecount1}) in the following way
\begin{eqnarray}\label{eq:posfacecount2}
p^{+,s}_{err}(k,m,n) & = & 2\sum_{l=m+2j+1,j\in \mN_0}^{n-1} \sum_{T^{(l)}\in \calF^{(l,+)},T^{(k-1)}\in T^{(l)}}\phiint(T^{(k-1)},T^{(l)})\phiext(T^{(l)},T^{(n-1)}) \nonumber \\
& = & 2\sum_{l=m+2j+1,j\in \mN_0}^{n-1} c^{(l,+)}_1 \phiint(T^{(k-1)},T^{(l)})\phiext(T^{(l)},T^{(n-1)}).
\end{eqnarray}
The key point in utilizing the projection of the standard simplexes is that for them the angles $\phiint(T^{(k-1)},T^{(l)})$, $\phiext(T^{(l)},T^{(n-1)})$ are known through the work of \cite{BorockyHenk} and \cite{Ruben} (see also, e.g. \cite{BetkeHenk,AS}), respectively. From \cite{BorockyHenk} we have for $\phiint(T^{(k)},T^{(l)})$
\begin{eqnarray}\label{eq:posfacecount3}
\phiint(T^{(k)},T^{(l)}) & = & \sqrt{\frac{l+1}{k+2}}\lp\frac{k+1}{k+2}\rp^{\frac{l-k-1}{2}}\lp\frac{k+2}{\pi}\rp^{\frac{l-k}{2}}
\frac{1}{\sqrt{\pi}}\int_{-\infty}^{\infty} \lp \frac{\sqrt{\pi}e^{-\frac{\lambda^2}{k+1}}\lp 1+i\erfi\lp\frac{\lambda}{\sqrt{k+1}}\rp\rp}{2\sqrt{k+1}}\rp^{l-k}e^{-\lambda^2}d\lambda \nonumber \\
& = &
\frac{\sqrt{l+1}}{2^{l-k}\sqrt{2\pi}}\int_{-\infty}^{\infty} e^{-\frac{\lambda^2(l+1)}{2}}\lp 1+i\erfi\lp\frac{\lambda}{\sqrt{2}}\rp\rp^{l-k}d\lambda.
\end{eqnarray}
From (\ref{eq:posfacecount3}) one then easily has
\begin{eqnarray}\label{eq:posfacecount3a}
\phiint(T^{(k-1)},T^{(l)}) =
\frac{\sqrt{l+1}}{2^{l-k+1}\sqrt{2\pi}}\int_{-\infty}^{\infty} e^{-\frac{\lambda^2(l+1)}{2}}\lp 1+i\erfi\lp\frac{\lambda}{\sqrt{2}}\rp\rp^{l-k+1}d\lambda.
\end{eqnarray}
Comparing (\ref{eq:posint1anal11}) and (\ref{eq:posfacecount3a}) we also observe that $\phiint(0,F^{(l,+)}_1)=\phiint(T^{(k-1)},T^{(l)})$ (this is not necessarily of an overwhelming interest here, but it will be a fairly useful observation in the following sections when we discuss the standard $\ell_1$). Also, following \cite{Ruben} we have for $\phiext(T^{(l)},T^{(n-1)})$
\begin{eqnarray}\label{eq:posfacecount4}
\phiext(T^{(l)},T^{(n-1)})=\sqrt{\frac{l+1}{\pi}}\int_{-\infty}^{\infty}e^{-(l+1)x^2}\lp\frac{1}{2}\erfc(-x)\rp^{n-l-1}dx.
\end{eqnarray}
A combination of (\ref{eq:posint1anal2}), (\ref{eq:posfacecount2}), (\ref{eq:posfacecount3a}), and (\ref{eq:posfacecount4}) is then enough to determine $p^{+,s}_{err}(k,m,n)$. Moreover, similarly to \cite{DonohoSigned}, one has
\begin{eqnarray}\label{eq:posfacecount5}
p^{+,s}_{err}(k,m,n)\leq p^{+}_{err}(k,m,n)\leq p^{+,s}_{err}(k,m,n+1).
\end{eqnarray}
In Table \ref{tab:l1regnonpersimprtab1} we show the results that one can obtain based on what we presented above. We also show corresponding simulated results and observe an excellent agreement between the simulated and theoretical results.

\begin{table}[h]
\caption{Simulated and theoretical results for $p^{+,s}_{err}(k,m,n)$ and $p^+_{err}(k,m,n)$}\vspace{.1in}
\hspace{-0in}\centering
\small{
\begin{tabular}{||c||c|c|c||}\hline\hline
$(k,m,n)$                & $ (3,5,8) $ & $ (3,4,6) $ & $ (4,5,8) $ \\ \hline \hline
$\#$ of failures         & $ 4091 $ & $ 6463 $ & $ 52850 $  \\ \hline
$\#$ of repetitions      & $ 32104 $ & $ 42252 $ & $ 154707 $  \\ \hline \hline
$p^{+,s}_{err}(k,m,n)$ -- \bl{\textbf{simulation}}/$p^{+,s}_{err}(k,m,n)$ -- \textbf{theory} & $ \bl{\mathbf{0.1274}}/\mathbf{0.1265} $ & $ \bl{\mathbf{0.1530}}/\mathbf{0.1539} $ & $ \bl{\mathbf{0.3416}}/\mathbf{0.3401} $  \\ \hline \hline
$\#$ of failures         & $ 6393 $ & $ 9941 $ & $  19666 $  \\ \hline
$\#$ of repetitions      & $ 42297 $ & $ 43096 $ & $ 49801 $  \\ \hline \hline
$p^+_{err}(k,m,n)$ -- \bl{\textbf{simulation}}/$p^+_{err}(k,m,n)$ -- \textbf{theory}    & $ \bl{\mathbf{0.1511}}/\mathbf{0.1528} $ & $ \bl{\mathbf{0.2307}}/\mathbf{0.2305} $ & $ \bl{\mathbf{0.3949}}/\mathbf{0.3937} $ \\ \hline \hline
$\#$ of failures         & $ 4813 $ & $ 12166 $ & $ 22175 $  \\ \hline
$\#$ of repetitions      & $ 23247 $ & $ 39421 $ & $ 46375 $  \\ \hline \hline
$p^{+,s}_{err}(k,m,n+1)$ -- \bl{\textbf{simulation}}/$p^{+,s}_{err}(k,m,n+1)$ -- \textbf{theory} & $ \bl{\mathbf{0.2070}}/\mathbf{0.2091} $ & $ \bl{\mathbf{0.3086}}/\mathbf{0.3077} $ & $ \bl{\mathbf{0.4782}}/\mathbf{0.4789} $  \\ \hline \hline
\end{tabular}
}
\label{tab:l1regnonpersimprtab1}
\end{table}

\section{Standard $\ell_1$}
\label{sec:l1}

In this section we look at the standard $\ell_1$ from (\ref{eq:l1}) and analyze its performance in a finite dimensional scenario. To facilitate reading we will try to follow fairly closely the strategy that we designed in Section \ref{sec:posl1}. Some results that trivially follow through the presentation of Section \ref{sec:posl1} we will just state without too much detailing; on the other hand, those that are substantially different we will try to emphasize.

We start things off by recalling on a result that is pretty much at the heart of everything that follows. For the concreteness of the exposition and without loss of generality we will continue to assume that the elements $\x_{k+1},\x_{k+2},\dots,\x_{n}$ of $\x$ are equal to zero and that the elements $\x_{1},\x_{2},\dots,\x_k$ have fixed signs (for the easiness of writing we will assume that $\x_{1},\x_{2},\dots,\x_k$ are nonnegative; differently from Section \ref{sec:posl1}, this assumption is only for the presentation of the analysis and is not a knowledge that can be used in the algorithms' design, i.e. it is not a priori known). The following is proven in \cite{StojnicCSetam09,StojnicICASSP09} and, as mentioned above, is among the key unsung heros that enabled running the entire machinery developed overthere.
\begin{theorem}(\cite{StojnicCSetam09,StojnicICASSP09} Nonzero elements of $\x$ have fixed signs and location)
Assume that an $m\times n$ systems matrix $A$ is given. Let $\x$
be a $k$ sparse vector. Also let $\x_{k+1}=\x_{k+2}=\dots=\x_{n}=0$. Let the signs of $\x_{1},\x_{2},\dots,\x_k$ be fixed, say all positive. Further, assume that $\y=A\x$ and that $\w$ is
a $n\times 1$ vector. If
\begin{equation}
(\forall \w\in \textbf{R}^n | A\w=0) \quad  -\sum_{i=1}^{k} \w_i<\sum_{i=k+1}^{n}|\w_{i}|,
\label{eq:thmcond1}
\end{equation}
then the solutions of (\ref{eq:l0gen}) and (\ref{eq:l1}) coincide. Moreover, if
\begin{equation}
(\exists \w\in \textbf{R}^n | A\w=0) \quad  -\sum_{i=1}^{k} \w_i\geq \sum_{i=k+1}^{n}|\w_{i}|,
\label{eq:thmcond2}
\end{equation}
then there is an $\x$ from the above set of $\x$'s with fixed location of nonzero elements such that the solution of (\ref{eq:l0gen}) is not the solution of (\ref{eq:l1}).
\label{thm:thmregposcond}
\end{theorem}
Similarly to how we introduced $C^+_w$ in (\ref{eq:posdefSw}), we set
\begin{equation}
C_w\triangleq\{\w\in \mR^n| \quad -\sum_{i=1}^k \w_i\geq \sum_{i=k+1}^{n}|\w_{i}|\}.\label{eq:defSw}
\end{equation}
Following (\ref{eq:posanal1}), the failing condition given in (\ref{eq:thmcond2}) can be utilized for performance characterization of (\ref{eq:l1}) in the following way
\begin{equation}
p_{err}(k,m,n)=P(\exists \w\in \textbf{R}^n | A\w=0, \quad \mbox{and} \quad -\sum_{i=1}^{k} \w_i\geq \sum_{i=k+1}^{n}|\w_{i}|)
=P(G^{(sub)}\cap C_w\neq \emptyset).\label{eq:anal1}
\end{equation}
$p_{err}(k,m,n)$ will denote the probability that the solution of (\ref{eq:l1}) is not the $k$-sparse solution of (\ref{eq:l0gen}) and $G^{(sub)}$ is, as earlier, a uniformly randomly chosen subspace from $G_{n,n-m}$, the Grassmanian of all $(n-m)$-dimensional subspaces of $\mR^n$. As $C^+_w$, $C_w$ is also a polyhedral cone and we write analogously to (\ref{eq:posanal2}) the following (again, a direct consequence of a remarkable result of Santalo \cite{Santalo} (see also, e.g. \cite{PMM,AS,BG}))
\begin{equation}
P(G^{(sub)}\cap C_w\neq \emptyset)=2\sum_{l=m+2j+1,j\in \mN_0}^{n} \sum_{F^{(l)}\in \calF^{(l,+)}}\phiint(0,F^{(l)})\phiext(F^{(l)},C_w),\label{eq:anal2}
\end{equation}
where $\calF^{(l)}$ is the set of all $l$-faces of $C_w$ ($\phiint(\cdot,\cdot)$ and $\phiext(\cdot,\cdot)$ are, as earlier, the internal and external angles, respectively). Connecting (\ref{eq:anal1}) and (\ref{eq:anal2}) we then also have
\begin{equation}
p_{err}(k,m,n)=P(G^{(sub)}\cap C_w\neq \emptyset)=2\sum_{l=m+2j+1,j\in \mN_0}^{n} \sum_{F^{(l)}\in \calF^{(l)}}\phiint(0,F^{(l)})\phiext(F^{(l)},C_w).\label{eq:anal3}
\end{equation}
(\ref{eq:anal3}) is of course an excellent conceptual way to characterize the performance of the standard $\ell_1$. However, as earlier, one again faces the same type of challenge, namely, to have (\ref{eq:anal3}) be of any practical use, one should be able to compute the angles $\phiint(\cdot,\cdot)$ and $\phiext(\cdot,\cdot)$. Below we present a strategy how it can be done. As in section \ref{sec:posl1}, the presentation will be split into two parts, the first one that relates to the internal angles and the second one that deals with the external angles.

\subsection{Internal angles}
\label{sec:intang}

In this section we analyze the internal angles appearing in (\ref{eq:anal3}), i.e. $\phiint(0,F^{(l)})$. We will distinguish between two cases: 1) $l<n$ and 2) $l=n$. For $l<n$, we have for the set of all $l$-faces $\calF^{(l)}_1$
\begin{multline}
\calF^{(l)}_1 \triangleq\{\w\in \mR^n| \quad -\sum_{i=1}^k \w_i= (\oness)^T\w_{I_r\setminus I^{(l)}_1},\diag(\oness)\w_{I_r\setminus I^{(l)}_1}\geq 0,\w_{I^{(l)}_1}=0,\\I^{(l)}_1\subset I_r,|I^{(l)}_1|=n-l-1,\oness\in\{-1,1\}^{l-k+1}\}.\label{eq:intanal2}
\end{multline}
The cardinality of set $\calF^{(l)}_1$ is easily given by
\begin{equation}
c^{(l)}_1\triangleq |\calF^{(l)}_1|=2^{l-k+1}\binom{n-k}{n-l-1}, l\in \{k-1,k,\dots,n-1\}.\label{eq:intanal4}
\end{equation}
(\ref{eq:anal3}) can then be rewritten in the following way
\begin{eqnarray}
p_{err}(k,m,n) & = & 2\sum_{l=m+2j+1,j\in \mN_0}^{n}  \sum_{F^{(l)}_1\in \calF^{(l)}_1}\phiint(0,F^{(l)}_1)\phiext(F^{(l)}_1,C_w)\nonumber \\
& = & 2 ( \sum_{l=m+2j+1,j\in \mN_0}^{n-1} c^{(l)}_1\phiint(0,F^{(l)}_1)\phiext(F^{(l)}_1,C_w) +  \phiint(0,C_w)\phiext(C_w,C_w)),
\label{eq:intanal6}
\end{eqnarray}
where due to symmetry $F^{(l)}_1$ is basically any of the elements from set $\calF^{(l)}_1$. For the concreteness we choose $I^{(l)}_1=\{l+2,l+3,\dots,n\}$ and $\oness=\1_{(l-k+1)\times 1}$ and consequently have
\begin{equation}
F^{(l)}_1 =\{\w\in \mR^n| \quad -\sum_{i=1}^k \w_i= \sum_{i=k+1}^{l+1}\w_{i},\w_{k+1:l+1}\geq 0,\w_{l+2:n}=0\}, l\in \{k-1,k,\dots,n-1\}=F^{(l,+)}_1.\label{eq:intanal7}
\end{equation}
A combination of (\ref{eq:posintanal7}), (\ref{eq:posint1anal11}), and (\ref{eq:intanal7}) gives
\begin{eqnarray}
\phiint(0,F^{(l)}_1)=\phiint(0,F^{(l,+)}_1)   =   \frac{\sqrt{l+1}}{2^{l-k+1}\sqrt{2\pi}}\int_{-\infty}^{\infty} \lp1-i\erfi\lp\frac{t}{\sqrt{2}}\rp\rp^{l-k+1} e^{-\frac{(l+1)t^2}{2}}dt.
\label{eq:int1anal11}
\end{eqnarray}
To compute $\phiint(0,C_w)$ we rely on (\ref{eq:posint2anal8}) and (\ref{eq:posint2anal12}). To that end we have
\begin{eqnarray}
\phiint(0,C_w) &  = & \frac{1}{(2\pi)^{\frac{l}{2}}}\int_{-\1_{1\times k}\w_{1:k}-\1_{1\times n-k}|\w_{k+1:n}|\geq 0}e^{-\frac{\w^T\w}{2}}d\w\nonumber \\
&  = & \frac{2^{n-k}}{(2\pi)^{\frac{n}{2}}}\int_{-\1_{1\times n}\w\geq 0,\w_{k+1:n}\geq 0}e^{-\frac{\w^T\w}{2}}d\w\nonumber \\
& = & \frac{1}{2\pi}\lim_{x\rightarrow \infty}\lim_{\epsilon\rightarrow 0_+}(\int_{-\infty}^{-\epsilon} \lp1-i\erfi\lp\frac{t}{\sqrt{2}}\rp\rp^{n-k} e^{-\frac{lt^2}{2}}\frac{(1-e^{-itx})}{it}dt\nonumber \\
& & +\int_{\epsilon}^{\infty} \lp1-i\erfi\lp\frac{t}{\sqrt{2}}\rp\rp^{n-k} e^{-\frac{lt^2}{2}}\frac{(1-e^{-itx})}{it}dt).\label{eq:int2anal8}
\end{eqnarray}

\subsection{External angles}
\label{sec:extang}

There are two types of the external angles that we need to compute $\phiext(F^{(l)}_1,C_w)$ and $\phiext(C_w,C_w)$. Trivially, one has $\phiext(C_w,C_w)=1$. Below we present a mechanism that can be used to handle $\phiext(F^{(l)}_1,C_w)$.

\subsubsection{Handling $\phiext(F^{(l)}_1,C_w)$}
\label{sec:ext1ang}

As mentioned earlier, the external angles at a given face represent the content/fraction of the subspace taken by the positive hull of the outward normals to the hyperplanes that meet at the face. For face $F^{(l)}_1$ we then have that the corresponding positive hull is given as
\begin{equation}\label{eq:ext1anal1}
phull^{(l)}_{ext,1}\triangleq  -pos\lp\begin{bmatrix}
                                      -\1_{(l+1)\times 1} \\
                                      -\1^{(s,1)}_{(n-l-1)\times 1}
                                    \end{bmatrix},\begin{bmatrix}
                                      -\1_{(l+1)\times 1} \\
                                      -\1^{(s,2)}_{(n-l-1)\times 1}
                                    \end{bmatrix},\dots,\begin{bmatrix}
                                      -\1_{(l+1)\times 1} \\
                                      -\1^{(s,2^{n-l-1})}_{(n-l-1)\times 1}
                                    \end{bmatrix}\rp,\1^{(s,i)}_{(n-l-1)\times 1}\in \{-1,1\}^{n-l-1},
\end{equation}
and $\1^{(s,i)}_{(n-l-1)\times 1}\neq \1^{(s,j)}_{(n-l-1)\times 1}$ for any $i\neq j$ and $1\leq i,j\leq 2^{n-l+1}$.
As in Section \ref{sec:posl1}, we will again rely on the ``Gaussian coordinates in an orthonormal basis" strategy and due to symmetry work with $-phull^{(l)}_{ext,1}$ instead of $phull^{(l)}_{ext,1}$. Moreover, the basis that we will work in will be the columns of the following matrix
\begin{equation}
B^{(l)}_{ext,1}=\begin{bmatrix}
e_{l+2} &  e_{l+3} & \dots & e_{n} & \begin{bmatrix}
                                       -\1_{(l+1)\times 1} \\
                                       \0_{(n-l-1)\times 1}
                                     \end{bmatrix}\frac{1}{\sqrt{l+1}}
          \end{bmatrix}.\label{eq:ext1anal2}
\end{equation}
Clearly, $-phull^{(l)}_{ext,1}$ is indeed located in the subspace spanned by the columns of $B^{(l)}_{ext,1}$. The tricky part of course is to accurately describe $-phull^{(l)}_{ext,1}$ through the coordinates corresponding to the basis vectors from $B^{(l)}_{ext,1}$. We will again denote the coordinates by a vector $\g$ ($\g\in \mR^{n-l}$) and look for a set of $\g$'s, say $G^{(l)}_{ext,1}$, such that
\begin{equation}
G^{(l)}_{ext,1} =\{\g\in \mR^{n-l}| B^{(l)}_{ext,1}\g\in -phull^{(l)}_{ext,1} \}.\label{eq:ext1anal3}
\end{equation}
Now let
\begin{equation}
D^{(l)}_{ext,1}=\begin{bmatrix}
\begin{bmatrix}
                                      -\1_{(l+1)\times 1} \\
                                      -\1^{(s,1)}_{(n-l-1)\times 1}
                                    \end{bmatrix} & \begin{bmatrix}
                                      -\1_{(l+1)\times 1} \\
                                      -\1^{(s,2)}_{(n-l-1)\times 1}
                                    \end{bmatrix} & \dots & \begin{bmatrix}
                                      -\1_{(l+1)\times 1} \\
                                      -\1^{(s,2^{n-l-1})}_{(n-l-1)\times 1}
                                    \end{bmatrix}\end{bmatrix},\label{eq:ext1anal4}
\end{equation}
and
\begin{eqnarray}
-phull^{(l)}_{ext,1} & = & \{D^{(l)}_{ext,1}\g^{(D)}|\g^{(D)}\geq 0,\g^{(D)}\in \mR^{2^{n-l-1}}\}\nonumber \\
& = & \{\begin{bmatrix}
\begin{bmatrix}
                                      -\1_{(l+1)\times 1} \\
                                      -\1^{(s,1)}_{(n-l-1)\times 1}
                                    \end{bmatrix} & \begin{bmatrix}
                                      -\1_{(l+1)\times 1} \\
                                      -\1^{(s,2)}_{(n-l-1)\times 1}
                                    \end{bmatrix} & \dots & \begin{bmatrix}
                                      -\1_{(l+1)\times 1} \\
                                      -\1^{(s,2^{n-l-1})}_{(n-l-1)\times 1}
                                    \end{bmatrix}\end{bmatrix}\g^{(D)}|\g^{(D)}\geq 0,\g^{(D)}\in \mR^{2^{n-l-1}}\}.\nonumber \\\label{eq:ext1anal5}
\end{eqnarray}
From (\ref{eq:ext1anal3}) and (\ref{eq:ext1anal5}) we have
\begin{equation}
G^{(l)}_{ext,1} =\{\g\in \mR^{n-l}| \exists \g^{(D)}\in\mR^{2^{n-l-1}}, \g^{(D)}\geq 0, \quad \mbox{and} \quad B^{(l)}_{ext,1}\g= D^{(l)}_{ext,1}\g^{(D)} \}.\label{eq:ext1anal6}
\end{equation}
First $(l+1)$ equations in $B^{(l)}_{ext,1}\g= D^{(l)}_{ext,1}\g^{(D)}$ imply that
\begin{equation}
\g_{n-l}\frac{1}{\sqrt{l+1}}=\sum_{i=1}^{2^{n-l-1}}\g^{(D)}_{i}\geq 0.\label{eq:ext1anal7}
\end{equation}
For $j\in\{2,3,\dots,n-l\}$, $(l+j)$-th row of $D^{(l)}_{ext,1}$ and $2^{n-l-2}$ has exactly $2^{n-l-2}$ ones and $2^{n-l-2}$ minus ones. Let the locations of $1$'s be $I^{(j,1)}$ and the locations of $-1$'s $I^{(j,-1)}$. Then the $j$-th equation in $B^{(l)}_{ext,1}\g= D_{ext,1}\g^{(D)}$ and (\ref{eq:ext1anal7}) also imply that
\begin{eqnarray}
& & \g_{j-1}=\sum_{i\in I^{j,1}}\g^{(D)}_{i}-\sum_{i\in I^{j,-1}}\g^{(D)}_{i}=2\sum_{i\in I^{(j,1)}}\g^{(D)}_{i}-\g_{n-l}\frac{1}{\sqrt{l+1}}\nonumber \\
\Longleftrightarrow & & \g_{j-1}+\g_{n-l}\frac{1}{\sqrt{l+1}}=2\sum_{i\in I^{(j,1)}}\g^{(D)}_{i}\geq 0 \nonumber \\
\Longleftrightarrow & & \g_{j-1}\geq -\g_{n-l}\frac{1}{\sqrt{l+1}}
\label{eq:ext1anal8}
\end{eqnarray}
Also,
\begin{eqnarray}
& & \g_{j-1}=\sum_{i\in I^{j,1}}\g^{(D)}_{i}-\sum_{i\in I^{j,-1}}\g^{(D)}_{i}=-2\sum_{i\in I^{(j,-1)}}\g^{(D)}_{i}+\g_{n-l}\frac{1}{\sqrt{l+1}}\nonumber \\
\Longleftrightarrow & & \g_{j-1}-\g_{n-l}\frac{1}{\sqrt{l+1}}=-2\sum_{i\in I^{(j,-1)}}\g^{(D)}_{i}\leq 0 \nonumber \\
\Longleftrightarrow & & \g_{j-1}\leq \g_{n-l}\frac{1}{\sqrt{l+1}}
\label{eq:ext1anal8a}
\end{eqnarray}
A combination of (\ref{eq:ext1anal6}), (\ref{eq:ext1anal7}), (\ref{eq:ext1anal8}), and (\ref{eq:ext1anal8a}) then gives
\begin{equation}
G^{(l)}_{ext,1} =\{\g\in \mR^{n-l}| \g_{n-l}\geq 0, -\g_{n-l}\frac{1}{\sqrt{l+1}}\leq \g_{j-1}\leq \g_{n-l}\frac{1}{\sqrt{l+1}}, j\in\{2,3,\dots,n-l\} \}.\label{eq:ext1anal9}
\end{equation}
(\ref{eq:ext1anal9}) is actually a complete characterization of $G^{(l)}_{ext,1}$, i.e. there are no other restrictions of $\g$. To see this we can assume that the bottom $n-l-1$ rows of $D^{(l)}_{ext,1}$ are ordered systematically. Say, $I^{(2,1)}=\{1,2,\dots,2^{n-l-2}\}$, $I^{(3,1)}=\{1,2,\dots,2^{n-l-3},2^{n-l-2}+1,2^{n-l-2}+2,\dots,32^{n-l-3}\}$, etc. Then it is clear that as $j$ grows there will always be enough degrees of freedom to set $\g^{(D)}$ so that (\ref{eq:ext1anal8}) and (\ref{eq:ext1anal8a}) hold.

 Since we determined  the content of $-phull^{(l)}_{ext,1}$ through the coordinates $\g$ of the basis $B^{(l)}_{ext,1}$ we switch to working with the independent standard normal, i.e. Gaussian, coordinates. To that end we finally have for $\phiext(F^{(l)}_1,C_w)$
\begin{eqnarray}
\phiext(F^{(l)}_1,C_w) &  = & \frac{1}{(2\pi)^{\frac{n-l}{2}}}\int_{\g\in G^{(l)}_{ext,1}}e^{-\frac{\g^T\g}{2}}d\g \nonumber \\
& = & \frac{1}{(2\pi)^{\frac{l}{2}}}\int_{\g_{n-l}\geq 0}e^{-\frac{\g_{n-l}^2}{2}}\lp \prod_{j=2}^{n-l}\frac{1}{(2\pi)^{\frac{l}{2}}}\int_{ -\g_{n-l}\frac{1}{\sqrt{l+1}}}^{\g_{n-l}\frac{1}{\sqrt{l+1}}}e^{-\frac{\g_{j-1}^2}{2}}d\g_{j-1}\rp d\g_{n-l}\nonumber \\
& = & \frac{1}{(2\pi)^{\frac{l}{2}}}\int_{\g_{n-l}\geq 0}e^{-\frac{\g_{n-l}^2}{2}}\lp \frac{1}{2}\erfc\lp\frac{-\g_{n-l}}{\sqrt{2}\sqrt{l+1}}\rp
-\frac{1}{2}\erfc\lp\frac{\g_{n-l}}{\sqrt{2}\sqrt{l+1}}\rp\rp^{n-l-1} d\g_{n-l}.\nonumber \\
\label{eq:ext1anal10}
\end{eqnarray}
What we presented above is then enough to determine $p_{err}(k,m,n)$. The obtained results are summarized in the following theorem.
\begin{theorem}(Exact $\ell_1$'s performance characterization -- finite dimensions)
Let $A$ be an $m\times n$ matrix in (\ref{eq:system})
with i.i.d. standard normal components (or, alternatively, with the null-space uniformly distributed in the Grassmanian). Let
the unknown $\x$ in (\ref{eq:system}) be $k$-sparse. Further, let the locations and signs of the nonzero components $\x$ be arbitrarily chosen but fixed.
Let $p_{err}(k,m,n)$ be the probability that the solutions of (\ref{eq:l0gen}) and (\ref{eq:l1}) do \emph{not} coincide. Then
\begin{eqnarray}
p_{err}(k,m,n) & = & 2 ( \sum_{l=m+2j+1,j\in \mN_0}^{n-1} c^{(l)}_1\phiint(0,F^{(l)}_1)\phiext(F^{(l)}_1,C_w)+\phiint(0,C_w)\phiext(C_w,C_w)),
\label{eq:finalthm1}
\end{eqnarray}
where $c^{(l)}_1$, $\phiint(0,F^{(l)}_1)$, $\phiint(0,C_w)$, and $\phiext(F^{(l)}_1,C_w)$ are as given in (\ref{eq:intanal4}), (\ref{eq:int1anal11}), (\ref{eq:int2anal8}), and (\ref{eq:ext1anal10}),  respectively.
\label{thm:finalperr}
\end{theorem}

\subsection{Simulations and theoretical results -- standard $\ell_1$}
\label{sec:thnumresults}

\begin{figure}[htb]
\begin{minipage}[b]{.5\linewidth}
\centering
\centerline{\epsfig{figure=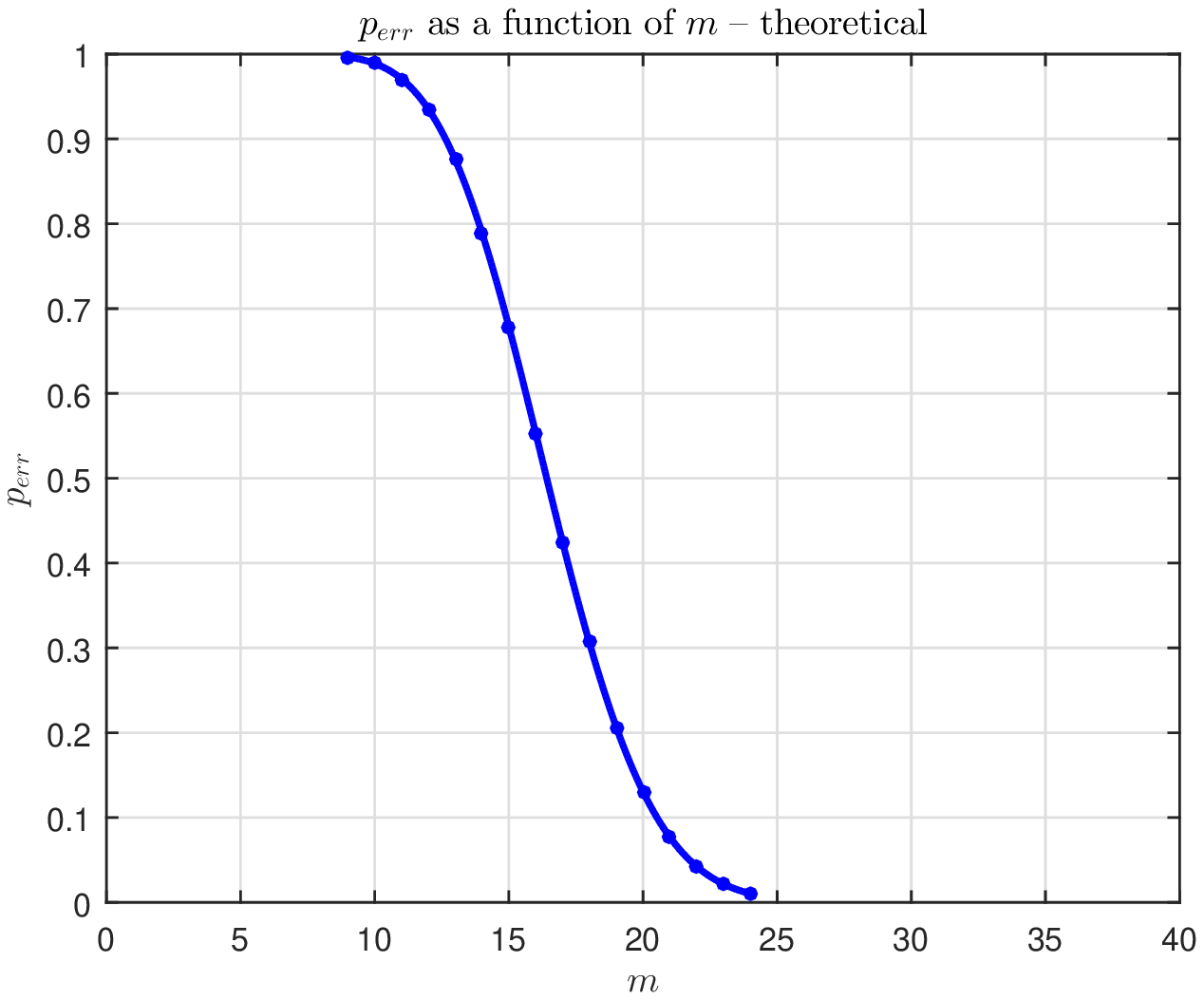,width=9cm,height=7cm}}
\end{minipage}
\begin{minipage}[b]{.5\linewidth}
\centering
\centerline{\epsfig{figure=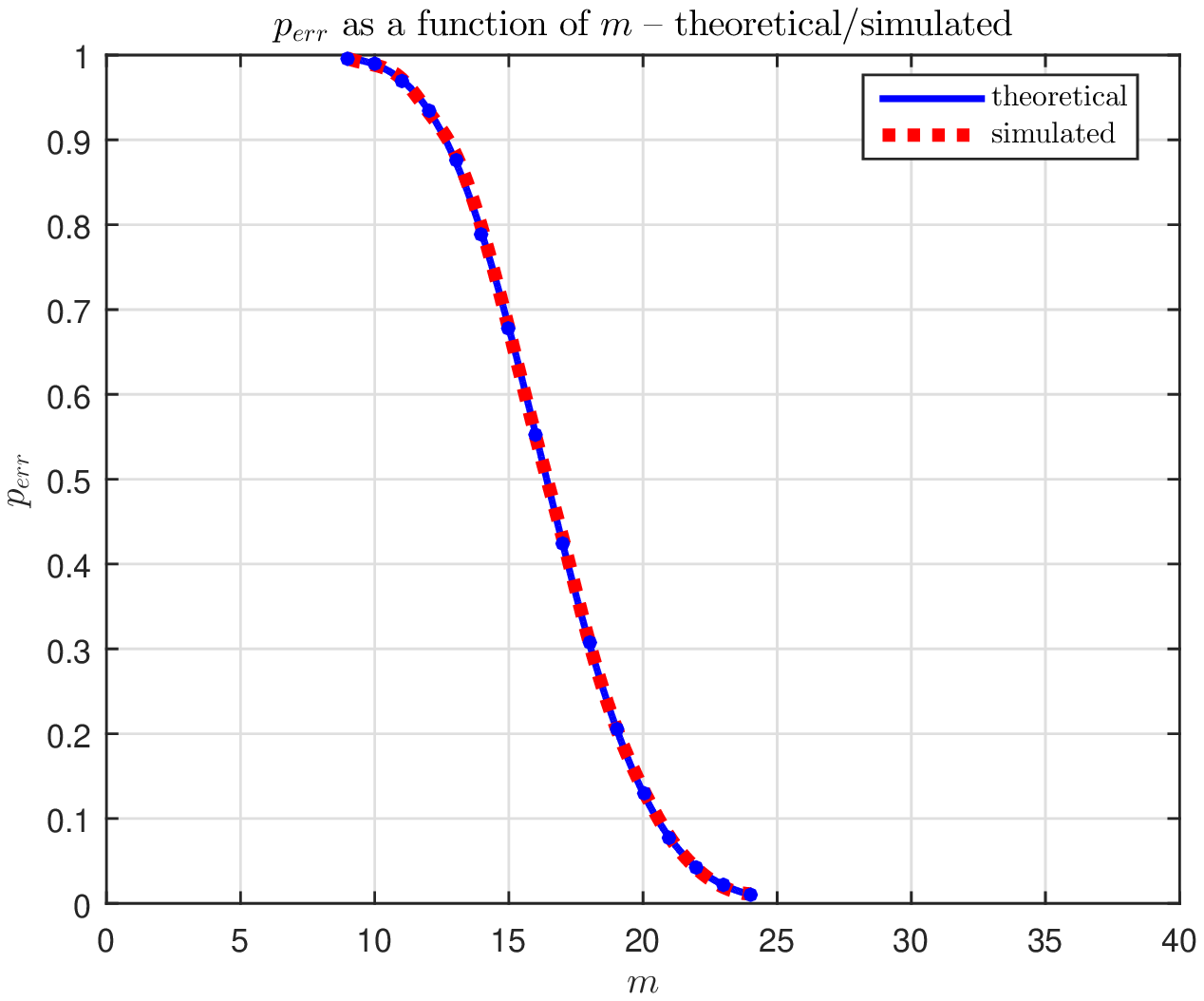,width=9cm,height=7cm}}
\end{minipage}
\caption{$p_{err}(k,m,n)$ as a function of $m$; left -- theory; right -- simulations}
\label{fig:l1regperr}
\end{figure}

In this section we present the concrete values for $p_{err}(k,m,n)$ that can be obtained through Theorem \ref{thm:finalperr} and numerical simulations. In Figure \ref{fig:l1regperr} we show both, the simulated and the theoretical values for $p_{err}(k,m,n)$ (the theoretical values are, of course, obtained based on Theorem \ref{thm:finalperr}). We fixed $k=6$ and $n=40$ and varied/increased $m$ so that $p_{err}(k,m,n)$ changes from one to zero. In addition to Figure \ref{fig:l1regperr} we also present Table \ref{tab:l1regperrtab1} where we show the numerical values for $p_{err}(k,m,n)$ (again, both, simulated and theoretical) obtained for several concrete values of triplets $(k,m,n)$ (we also show the number of numerical experiments that were run as well as the number of them that did not result in having the solution of (\ref{eq:l1}) match the $k$-sparse solution of (\ref{eq:l0gen})). As in Section \ref{sec:posl1}, we observe an excellent agreement between the simulated and theoretical results.

\begin{table}[h]
\caption{Simulated and theoretical results for $p_{err}(k,m,n)$; $k=6$, $n=40$}\vspace{.1in}
\hspace{-0in}\centering
\begin{tabular}{||c||c|c|c|c|c|c||}\hline\hline
$m$                    & $ 14 $ & $ 15 $ & $ 16 $ & $ 17 $ & $ 18 $ & $ 19 $ \\ \hline \hline
$\#$ of failures       & $ 11343 $ & $ 6777 $ & $ 2820 $ & $ 3246 $ & $ 2509 $ & $ 1166 $ \\ \hline
$\#$ of repetitions    & $ 14262 $ & $ 10011 $ & $ 5125 $ & $ 7523 $ & $ 8156 $ & $ 5625 $ \\ \hline \hline
$p_{err}$ -- simulation& $ \bl{\mathbf{0.7953}} $ & $ \bl{\mathbf{0.6770}} $ & $ \bl{\mathbf{0.5502}} $ & $ \bl{\mathbf{0.4315}} $ & $ \bl{\mathbf{0.3076}} $ & $ \bl{\mathbf{0.2073}} $ \\ \hline \hline
$p_{err}$ -- theory    & $ \mathbf{0.7890} $ & $ \mathbf{0.6786} $ & $ \mathbf{0.5530} $ & $ \mathbf{0.4248} $ & $ \mathbf{0.3064} $ & $ \mathbf{0.2069} $ \\ \hline \hline
\end{tabular}
\label{tab:l1regperrtab1}
\end{table}

\subsection{Asymptotics}
\label{sec:asym}

In this section we look at what happens in an asymptotic scenario. As in Section \ref{sec:posl1}, we will assume that he systems dimensions grow larger in a linearly proportional fashion. Similarly to the positive $\ell_1$, the standard $\ell_1$ also exhibits the so-called phase-transition (PT) phenomenon. That means that for a given $\alpha\in(0,1)$ there will be a $\beta_w\in(0,\alpha)$ such that for $\beta< \beta_w$ the $k$-sparse solution of (\ref{eq:l0gen}) is the solution of (\ref{eq:l1}) with overwhelming probability. On the other hand, for $\beta\geq \beta_w$ there will be a $k$-sparse $\x$, from a given set of $\x$'s with randomly chosen but fixed locations and signs, that is the solution of (\ref{eq:l0gen}) and with overwhelming probability is not the solution of (\ref{eq:l1}). As usual, a full characterization of the PT phenomenon goes through the determination of the phase-transition curve (PT curve) in $(\alpha,\beta)$ plane, i.e. through the determination of the so-called breaking points where, roughly speaking, (\ref{eq:l1}) solves or does not solve (\ref{eq:l0gen}). As in the case of the positive $\ell_1$, one can also conduct the large deviations principle (LDP) type of analysis to obtain a way more complete picture that describes the behavior of (\ref{eq:l1}) in the entire transition zone around the breaking points. This typically goes through the determination of the so-called rate of decay of the probability of success/failure, i.e. through the determination of $I_{ldp}(\alpha,\beta)=\lim_{n\rightarrow\infty}\frac{\log(p_{err}(k,m,n))}{n}$. Settling the LDP behavior of (\ref{eq:l1}) assumes computing $I_{ldp}(\alpha,\beta)$ for any $(\alpha,\beta)\in(0,1)\times (0,1)$ for which the original system has unique solution.

The analysis that we presented in earlier sections can be used as a starting point on a path to determine the $\ell_1$'s LDP behavior and consequently its PT curve as well. Since settling the PT phenomenon alone is much easier than handling the entire LDP it can also be done in separate ways that don't rely on the LDP results. As is of course well known by now, a full characterization of the standard $\ell_1$'s PT behavior has been obtained through the work of \cite{DonohoPol,DonohoUnsigned,StojnicCSetam09,StojnicUpper10}. Settling LDP, as a harder problem, required a bit extra effort and was finally done in \cite{Stojnicl1RegPosasymldp}. The LDP analysis presented in \cite{Stojnicl1RegPosasymldp} is somewhat involved and a detailed discussion about it goes way beyond what we study here. Nonetheless, we below briefly sketch how one can bridge between the above analysis and what is presented in \cite{Stojnicl1RegPosasymldp}. In other words, we below show what problems a further analysis of the above $p_{err}(k,m,n)$ characterizations boils down to in an infinite dimensional setting (precisely these problems are, among other things, what is further discussed in \cite{Stojnicl1RegPosasymldp}).

Assuming the linear asymptotic regime (i.e. $k=\beta n, m=\alpha n$, and $l=\rho n$, where $\beta$, $\alpha$, and $\rho$ are fixed constants independent of $n$) as $n\rightarrow\infty$, from (\ref{eq:finalthm1}) we have
\begin{equation}
\lim_{n\rightarrow \infty}\frac{\log(p_{err}(k,m,n))}{n}  =  \max\{\max_{\rho\geq \alpha}\lim_{n\rightarrow \infty} \frac{\log(\zeta^{(\infty)}_1)}{n},\lim_{n\rightarrow \infty}\frac{\log(\zeta^{(\infty)}_2)}{n}\},
\label{eq:asym1}
\end{equation}
where
\begin{eqnarray} \label{eq:asym2}
  \lim_{n\rightarrow \infty} \frac{\log(\zeta^{(\infty)}_1)}{n} &=&
  \lim_{n\rightarrow \infty} \lp\frac{\log(c^{(l)}_{1})}{n}+
  \frac{\log(\phiint(0,F^{(l)}_1))}{n}+\frac{\log(\phiext(F^{(l)}_1,C_w))}{n}\rp \nonumber \\
  \lim_{n\rightarrow \infty} \frac{\log(\zeta^{(\infty)}_2)}{n} &=&
  \lim_{n\rightarrow \infty} \lp
  \frac{\log(\phiint(0,C_w))}{n}+\frac{\log(\phiext(C_w,C_w))}{n}\rp.
\end{eqnarray}
From (\ref{eq:intanal4}) we recall that
\begin{eqnarray}
\lim_{n\rightarrow \infty} \frac{\log(c^{(l)}_{1})}{n}=\lim_{n\rightarrow \infty}\frac{\log\lp 2^{l-k+1}\binom{n-k}{n-l-1}\rp}{n}=-(1-\beta)H\lp\frac{1-\rho}{1-\beta}\rp +(\rho-\beta)\log(2).\label{eq:asym3}
\end{eqnarray}
Let us also recall that
\begin{eqnarray}
\phiint(0,F^{(n-1)}_1) &  = &
\frac{\sqrt{n}}{(2\pi)^{\frac{n-1}{2}}}\int_{-\1_{1\times n-1}\w_{1:n-1}\geq 0,\w_{k+1:n-1}\geq 0}e^{-\frac{\w_{1:n-1}^T\w_{1:n-1}+\w_{1:n-1}^T( -\1_{1\times n-1})^T
             (-\1_{1\times n-1})\w_{1:n-1}}{2}}d\w_{1:n-1}.\nonumber \\
\label{eq:asym4a}
\end{eqnarray}
Now, although $F^{(n)}_1$ is not a face of $C_w$, let us define a mathematical object $\phiint(0,F^{(n)}_1)$ in the following way (essentially a mathematically analogous way to (\ref{eq:asym4a}))
\begin{eqnarray}
\phiint(0,F^{(n)}_1) &  = &
\frac{\sqrt{n+1}}{(2\pi)^{\frac{n}{2}}}\int_{-\1_{1\times n}\w_{1:n}\geq 0,\w_{k+1:n}\geq 0}e^{-\frac{\w_{1:n}^T\w_{1:n}+\w_{1:n}^T( -\1_{1\times n})^T
             (-\1_{1\times n})\w_{1:n}}{2}}d\w_{1:n}.\nonumber \\
\label{eq:asym4b}
\end{eqnarray}
Then we have
\begin{eqnarray}
\lim_{n\rightarrow \infty} \frac{\log(\phiint(0,F^{(n-1)}_1)}{n} &  = & \lim_{n\rightarrow \infty} \frac{\log(\phiint(0,F^{(n)}_1)}{n}\nonumber \\
&  \geq  & \lim_{n\rightarrow \infty} \frac{\log\lp\frac{1}{(2\pi)^{\frac{n}{2}}}\int_{-\1_{1\times n}\w_{1:n}\geq 0,\w_{k+1:n}\geq 0}e^{-\frac{\w_{1:n}^T\w_{1:n}}{2}}d\w_{1:n}\rp}{n}\nonumber \\
&  =  & \lim_{n\rightarrow \infty} \frac{\log\lp\frac{1}{2^{n-k}}\phiint(0,C_w)\rp}{n},\nonumber \\\label{eq:asym4c}
\end{eqnarray}
where the inequality (basically an equality) follows by a combination of (\ref{eq:posasym4h}), (\ref{eq:posasym4i}), and (\ref{eq:int1anal11}). From (\ref{eq:asym1}), (\ref{eq:asym2}), (\ref{eq:asym3}), and (\ref{eq:asym4c}) we obtain
\begin{equation}
\lim_{n\rightarrow \infty}\frac{\log(p_{err}(k,m,n))}{n}  =  \max_{\rho\geq \alpha}\lim_{n\rightarrow \infty} \frac{\log(\zeta^{(\infty)}_1)}{n}.
\label{eq:asym5}
\end{equation}
A further combination of (\ref{eq:posasym4h}), (\ref{eq:posasym4i}), (\ref{eq:posasym7}), (\ref{eq:int1anal11}), (\ref{eq:ext1anal10}), (\ref{eq:asym2}), and (\ref{eq:asym3}) transforms (\ref{eq:asym5}) into the following
\begin{eqnarray}
\lim_{n\rightarrow \infty}\frac{\log(p_{err}(k,m,n))}{n}
 & = & \max_{\rho\geq \alpha}\lim_{n\rightarrow \infty} \frac{\log(\zeta^{(\infty)}_1)}{n}\nonumber \\
 & = & \max_{\rho\geq \alpha} (-(1-\beta)H\lp\frac{1-\rho}{1-\beta}\rp+(\rho-\beta)\log(2) \nonumber \\
 & & + \min_{\mu_y\geq 0} \lp(\rho-\beta) \log(\erfc(\mu_y))+\rho\mu_y^2\rp-(\rho-\beta)\log(2)\nonumber \\.
& & + \max_{\g_{n-l}\geq 0} \lp -\rho \g_{n-l}^2+(1-\rho)\log\lp \frac{1}{2}\erfc(-\g_{n-l})-\frac{1}{2}\erfc(\g_{n-l})\rp\rp).
\label{eq:asym12}
\end{eqnarray}
For a given $\beta$, let $\alpha_w$ be the $\alpha$ that gives $\lim_{n\rightarrow \infty}\frac{\log(p_{err}(k,m,n))}{n} =0$ (as mentioned after (\ref{eq:posasym12}) such an $\alpha$ always exists). One can then also repeat the rest of the reasoning after (\ref{eq:posasym12}) to conclude that if $\alpha\geq \alpha_w$ then $\rho=\alpha$ is optimal in (\ref{eq:asym12}) and  the optimization over $\rho$ in (\ref{eq:asym12}) can be removed.
Also, the same reasoning applies if $\alpha\leq \alpha_w$ and the optimal $\rho$ will again trivially be equal to $\alpha$. The only difference is that when $\alpha\leq \alpha_w$ one focuses on the complementary version of (\ref{eq:anal2})
\begin{equation}
P(G^{(sub)}\cap C_w\neq \emptyset)=1-2\sum_{l=m-2j-1,j\in \mN_0,l\geq k-1} \sum_{F^{(l)}\in \calF^{(l)}}\phiint(0,F^{(l)})\phiext(F^{(l)},C_w)=1-p_{cor},\label{eq:asym13}
\end{equation}
where $p_{cor}$ is the probability of being correct, i.e. the probability that the solution of (\ref{eq:l1}) is the $k$-sparse solution of (\ref{eq:l0gen}) and its decay rate is given by (\ref{eq:asym12}) with $\rho\leq \alpha$. (\ref{eq:asym12}) is then sufficient to fully determine numerically PT and LDP curves of the standard $\ell_1$. As was the case for the positive $\ell_1$, \cite{Stojnicl1RegPosasymldp}, however, goes way beyond that and determines explicit solutions to (\ref{eq:asym12}).

\subsection{Face counting and projection survival}
\label{sec:facecountprsur}

The connection to the counting of the faces and survival of their projections that we discussed in Section \ref{sec:posfacecountprsur} can be established for the standard $\ell_1$ as well. As mentioned in Section \ref{sec:posfacecountprsur}, in \cite{DonohoUnsigned} Donoho established roughly the following: if and only if a given $(k-1)$-face of the regular standard $n$-dimensional crosspolytope $C^{(n)}$ (that essentially corresponds to  a set of scaled $k$-sparse $\x$ in (\ref{eq:system})) so to say ``survives" the projection by $A$, i.e. remains a face of the projected crosspolytope $C^{(n)}$, $AC^{(n)}$, then the solution of (\ref{eq:l1}) is the $k$-sparse solution of (\ref{eq:l0gen}) (which is located exactly on the given $(k-1)$-face of $C^{(n)}$). In \cite{DonohoPol} Donoho then switched to the problem of counting the faces of $AC^{(n)}$. He relied on a relation due to Affentrager and Schneider that connects the numbers of faces of various geometric bodies and their random projections. This relation at its core uses the probability that a particular face survives the projection (i.e. remains a face after the projection) which was introduced in \cite{Grunbaum68} as the Grunbaum's Grassman angle. If one is interested in characterizing the performance of the standard $\ell_1$ and not necessarily in the count of the projected faces, Donoho's approach of relying on \cite{AS} can be circumvented and instead one can simply observe that
\begin{equation}\label{eq:facecount0}
p^{FC}_{err}(k,m,n)  =\gamma_{GB}(F^{(k-1)},C^{(n)}),
\end{equation}
where $\gamma_{GB}(F^{(k-1)},C^{(n)})$ is the above mentioned Grunbaum's Grassman angle and $F^{(k-1)}$ is a $(k-1)$-face of $C^{(n)}$. A very nice thing about the Grunbaum's Grassman angle is its
a characterization through \cite{Santalo} and the angle-sum relation of \cite{PMM} which boils down to the following
\begin{equation}\label{eq:facecount1}
p^{FC}_{err}(k,m,n)  =\gamma_{GB}(F^{(k-1)},C^{(n)})  =2\sum_{l=m+2j+1,j\in \mN_0}^{n} \sum_{F^{(l)}\in \calF^{(l)}_{C^{(n)}},F^{(k-1)}\in F^{(l)}}\phiint(F^{(k-1)},F^{(l)})\phiext(F^{(l)},C^{(n)}),
\end{equation}
where $\calF^{(l)}_{C^{(n)}}$ is the set of all $l$-faces of $C^{(n)}$. As is well known, for $l<n$, $F^{(l)}=T^{(n)}$, and for $l=n$ we have $F^{(l)}=C^{(n)}$. Also, as observed in \cite{DonohoPol,BorockyHenk}, if $l<n$, the number of $l$-faces that contain the given $(k-1)$-face is equal to $c^{(l)}_1$ and one can rewrite (\ref{eq:facecount1}) in the following way
\begin{eqnarray}\label{eq:facecount2}
p^{FC}_{err}(k,m,n) & = & 2\sum_{l=m+2j+1,j\in \mN_0}^{n} \sum_{T^{(l)}\in \calF^{(l)}_{C^{(n)}},T^{(k-1)}\in T^{(l)}}\phiint(T^{(k-1)},T^{(l)})\phiext(T^{(l)},C^{(n)})   \nonumber \\
& = & 2(\sum_{l=m+2j+1,j\in \mN_0}^{n-1} c^{(l)}_1 \phiint(T^{(k-1)},T^{(l)})\phiext(T^{(l)},C^{(n)}) +\phiint(T^{(k-1)},C^{(n)})\phiext(C^{(n)},C^{(n)})).\nonumber \\
\end{eqnarray}
We also add that, trivially, for $l=n$ there is only one face (basically $C^{(n)}$) that contains any given $(k-1)$-face. The key point in counting the numbers of the projected standard crosspolytope faces in \cite{BorockyHenk} and in asymptotically characterizing the performance of the $\ell_1$ in \cite{DonohoPol,DonohoUnsigned} is that for the standard crosspolytopes all the relevant angles needed here are known. Namely, $\phiint(T^{(k-1)},T^{(l)})$ and $\phiext(T^{(l)},C^{(n)})$ are known through the work of \cite{BetkeHenk} (see also, e.g. \cite{BorockyHenk,AS}). Strictly speaking, in a nonasymptotic regime when $m+1$ and $n$ are both even or both odd, (\ref{eq:facecount2}) requires one to know $\phiint(T^{(k-1)},C^{(n)})$ as well; however one can instead of (\ref{eq:facecount2}) use then the characterization of the Grunbaum's Grassman angle through the bottom portion of the McMullen's angle sum, e.g.
\begin{equation}\label{eq:facecount3}
p^{FC}_{err}(k,m,n)  = 1- 2\sum_{l=m-2j-1,j\in \mN_0}^{n} \sum_{T^{(l)}\in \calF^{(l)}_{C^{(n)}},T^{(k-1)}\in T^{(l)}}\phiint(T^{(k-1)},T^{(l)})\phiext(T^{(l)},C^{(n)}),
\end{equation}
so that knowing $\phiint(T^{(k-1)},T^{(l)})$ and $\phiext(T^{(l)},C^{(n)})$ is indeed sufficient. In \cite{BetkeHenk}, (\ref{eq:posfacecount3a}) was established
\begin{eqnarray}\label{eq:facecount3a}
\phiint(T^{(k-1)},T^{(l)}) =
\frac{\sqrt{l+1}}{2^{l-k+1}\sqrt{2\pi}}\int_{-\infty}^{\infty} e^{-\frac{\lambda^2(l+1)}{2}}\lp 1+i\erfi\lp\frac{\lambda}{\sqrt{2}}\rp\rp^{l-k+1}d\lambda,
\end{eqnarray}
and, as observed earlier, $\phiint(0,F^{(l)}_1)=\phiint(T^{(k-1)},T^{(l)})$. Also, following \cite{BetkeHenk} we have for $\phiext(T^{(l)},C^{(n)})$
\begin{eqnarray}\label{eq:facecount4}
\phiext(T^{(l)},C^{(n)})=\sqrt{\frac{l+1}{\pi}}\int_{0}^{\infty}e^{-(l+1)x^2}\lp\erf(x)\rp^{n-l-1}dx.
\end{eqnarray}
Comparing (\ref{eq:ext1anal10}) and (\ref{eq:facecount4}) it is not that hard to see that $\phiext(F^{(l)}_1,C_w)=\phiext(T^{(l)},C^{(n)})$. Although it is not necessarily needed, we mention that one can compute $\phiint(T^{(k-1)},C^{(n)})$ explicitly as well; it does turn out that $\phiint(T^{(k-1)},C^{(n)})=\phiint(0,C_w)$. A combination of (\ref{eq:facecount2}), (\ref{eq:facecount3}), (\ref{eq:facecount3a}), and (\ref{eq:facecount4}) is then enough to determine $p^{FC}_{err}(k,m,n)$. As it must be, one has
\begin{eqnarray}\label{eq:facecount5}
p^{FC}_{err}(k,m,n)= p_{err}(k,m,n).
\end{eqnarray}
Of course, if one is solely interested in the performance of the standard $\ell_1$, the results from \cite{DonohoPol} (and ultimately the angle characterizations from \cite{BorockyHenk,BetkeHenk} and their connections to the Grunbaum's Grassman angle from \cite{AS,PMM,Grunbaum68,Santalo}) could have been used. In other words, there would be no need to adapt the positive $\ell_1$ results from Section \ref{sec:posl1}. However, we found the methods from Section \ref{sec:posl1} generically very useful when studying a large class of other problems and in this introductory paper wanted to showcase how easily they can handle the standard $\ell_1$ case as well.

\section{Conclusion}
\label{sec:conc}

In this paper we studied random linear systems with sparse solutions. In particular, we provided a performance characterization of the well known $\ell_1$ heuristic. Instead of typical large (theoretically infinite) dimensional settings we here focused on the finite dimensional scenarios. We designed several novel high-dimensional geometry techniques that turned out to be of great help. We fully characterized the so-called positive/nonnegative $\ell_1$ (applicable in scenarios where one knows in advance that the unknown vectors are not only sparse but also nonnegative) and then showed how one can adapt such a characterization so that it fits the standard $\ell_1$ as well. For both of these problems we also showed how to adapt the obtained results as one moves towards the asymptotic regime. We obtained explicit formulations of the problems one would need to solve to fully characterize such a behavior as well. Of course, a typical way of handling the asymptotic regime goes through the understanding of the so-called phase-transition phenomenon. In our companion paper \cite{Stojnicl1RegPosasymldp}, we among other things rely on the asymptotic formulations provided here to go way beyond the typical phase-transition phenomenon and provide a full characterization of much broader large deviations concepts.

Also, for both, the positive and the standard $\ell_1$, we revisited their connection with the face counting problems in high-dimensional geometry and provided several interesting observations in that direction as well. Finally, we presented a solid collection of results obtained through numerical simulations for all problems considered here and observed an excellent level of agreement between them on the one side and what the theory predicts on the other.

As expected, a decent level of simplicity and elegance in the used arguments leave a tone of opportunity to continue further and consider various other aspects/extensions of the problems considered here. This usually assumes fairly routine adjustments of the techniques introduced here. For a couple of interesting problems we will in a few companion papers present how these adjustments can be done and what results they eventually produce.

In addition to those routine types of further considerations, we would like to single out one that, in our view, stands out. Namely, we designed a large number of novel techniques that turned out to be very powerful in handling the asymptotic behavior of many hard problems. For them the Gaussian distribution turns out to be the common denominator, i.e. pretty much any distribution that can be pushed through the central limit theorem in asymptotic scenario typically produces the same results as does the Gaussian. In the final dimensional setting considered here that does not seem to be the case. We believe that it is an extraordinary challenge to provide an analogue of the results presented here for types of statistics different from the Gaussian (of course, of the real interests would be those where the final results substantially differ from the Gaussian ones presented here).

\begin{singlespace}
\bibliographystyle{plain}
\bibliography{l1regposfinn1Refs}
\end{singlespace}

\end{document}